\documentclass[11pt,a4paper]{article}

\usepackage{mathtools}
\usepackage{authblk} % author
\usepackage{algpseudocode,algorithmicx,algorithm}

\usepackage{mathrsfs}
\usepackage{latexsym,bm}
\usepackage{amsmath,amsfonts,amsmath,amssymb,amsthm}
\usepackage{extarrows}

\usepackage{graphicx,subfigure,epstopdf,float}
\usepackage{enumerate,cases,multirow}
\usepackage{caption}

\usepackage{longtable,colortbl,arydshln,threeparttable}
\definecolor{mygray}{gray}{.9}

\usepackage{indentfirst}
\usepackage[top=25mm,bottom=20mm,left=25mm,right=20mm]{geometry}
\usepackage{xcolor}
\usepackage{cite}

\usepackage{listings}

\usepackage{makeidx}        % generate an index, automatically
\usepackage{booktabs}
\usepackage[bookmarksnumbered,colorlinks,citecolor=red,linkcolor=red,hyperindex,linktocpage=true]{hyperref}

\newcommand{\ket}[1]{| #1 \rangle} % |u>
\newcommand{\bra}[1]{\langle #1 |} % <u|

\newcommand{\bb}{\boldsymbol}

\def \d {\mathrm{d}}
\def \e {\mathrm{e}}
\def \i {\mathrm{i}}

\newcounter{parentalgorithm}

\makeatother

\newtheorem{theorem}{Theorem}[section]
\newtheorem{lemma}{Lemma}[section]

\newtheorem{definition}{Definition}[section]

\theoremstyle{remark}
\newtheorem{remark}{\bf Remark}[section]

\numberwithin{equation}{section}

\usepackage{lineno}

\begin{document}
%\begin{CJK*}{GBK}{song}

%\pagewiselinenumbers % Renumber by page
%\linenumbers
%%\switchlinenumbers

\title{Time complexity analysis of  quantum difference methods for linear high dimensional and multiscale  partial differential equations}
\author[1]{Shi Jin \thanks{shijin-m@sjtu.edu.cn}}
\author[2,3,4]{Nana Liu \thanks{nana.liu@quantumlah.org}}
\author[1]{Yue Yu \footnote{Corresponding author.} \thanks{terenceyuyue@sjtu.edu.cn}}
\affil[1]{School of Mathematical Sciences, Institute of Natural Sciences, MOE-LSC, Shanghai Jiao Tong University, Shanghai, 200240, P. R. China.}
\affil[2]{Institute of Natural Sciences, Shanghai Jiao Tong University, Shanghai 200240, China.}
\affil[3]{Ministry of Education, Key Laboratory in Scientific and Engineering Computing, Shanghai Jiao Tong University,
Shanghai 200240, China}
\affil[4]{University of Michigan-Shanghai Jiao Tong University Joint Institute, Shanghai 200240, China}
\date{}
\maketitle

\begin{abstract}
We investigate time complexities of finite difference methods for solving the high-dimensional linear heat equation, the high-dimensional linear  hyperbolic equation and the multiscale hyperbolic heat system with quantum algorithms (hence referred to as the ``quantum difference methods'').  For the heat and linear hyperbolic equations we study the impact of explicit and implicit time discretizations on quantum advantages over the classical difference method. For the multiscale problem, we find  the time complexity of both the classical treatment and quantum treatment for the explicit scheme scales as $\mathcal{O}(1/\varepsilon)$, where $\varepsilon$ is the scaling parameter,  while the scaling for the multiscale  Asymptotic-Preserving (AP) schemes does not depend on $\varepsilon$. This indicates that it is still of great importance to develop AP schemes for multiscale problems in quantum computing.
\end{abstract}


\tableofcontents

\section{Introduction}

High-dimensional problems and multiscale problems have been two challenges in scientific computing. For high-dimensional partial differential equations (PDEs), such as the $N-$body Schr\"odinger equation in quantum mechanics, the Boltzmann equation in kinetic theory, and PDEs with (high-dimensional) uncertainties \cite{UQ-CFD, JinPareschi-book}, classical numerical methods suffer from the so-called ``curse of dimensionality'' since the computational cost often increases exponentially with the dimension of the problem, which is undesirable for simulations on classical computers. To overcome this bottleneck, Monte-Carlo methods, sparse grid methods, and/or mean-field approximations are often used computational or mathematical tools.
In addition, in these problems there are often  multiple  spatial and temporal scales that pose further challenges for numerical computations. This is because a naive numerical discretization of these equations requires the mesh sizes and time steps smaller than the scaling parameters, which is prohibitively expensive in applications \cite{weinan2011principles}.
Among various multiscale methods, the Asymptotic-Preserving (AP) scheme, which preserves the asymptotic transition from the micro models to the macro ones at the discrete level,  has the merit of using one solver that works across scales naturally, thus has been widely used for multiscale hyperbolic and kinetic problems  \cite{JinReview-2012, JinActa}.

In this paper, we consider the application of quantum linear systems algorithms (QLSA) to these two types of problems, namely high dimensional and multiscale PDEs. Quantum computing is a rapidly growing computational paradigm  that has attracted significant attention due to the discovery of quantum algorithms that have exponential acceleration over the best-known classical methods \cite{Nielsen2010,Lipton2010,Liu-Thompson-Weedbrook-2016,Liu-Demarie-Tan-2019,Deutsch-1992,Shor-1997,HHL-2009}. In \cite{HHL-2009}, Harrow, Hassidim and Lloyd propose a quantum algorithm (HHL algorithm) for solving linear system of equations and proved that the algorithm can provide exponential speedup over the classical conjugate gradient (CG) method in terms of the matrix size. Later Cao et al. \cite{CaoCircuit-2012} present an efficient and generic quantum circuit design for implementing the algorithm. The exponential speedup of the HHL algorithm is expected to break the curse of dimensionality, which has attracted research to apply the method to solve  linear systems resulting from  classical numerical discretizations of both ordinary and partial differential equations.
For example, Berry in \cite{Berry-2014} first discretize the first-order linear ordinary differential equation (ODE) by using the linear multi-step method and then apply the HHL algorithm to solve the resulting linear system of equations. The HHL algorithm requires $\mathcal{O}(1/\delta)$ uses of a unitary operation to estimate its eigenvalues to precision $\delta$. To circumvent the limitations of phase estimation, Childs et al. propose a new algorithm to exponentially improve the dependence on the precision parameter in \cite{Childs-2017} and apply the quantum method to develop a finite difference algorithm for the Poisson equation and a spectral algorithm for more general second-order elliptic equations in \cite{Childs-2021}, respectively.
However, even if this improved QLSA is used to implement the quantum difference method in \cite{Berry-2014}, the overall complexity is still $\text{poly}(1/\delta)$ since the multi-step method itself is a significant source of error \cite{Berry-Childs-Ostrander-2017}.
To circumvent this limitation, Berry, Childs and Ostrander et al. present a quantum algorithm for linear differential equations
with complexity $\text{polylog}(1/\delta)$, where a truncation of the Taylor series of the propagator for the differential equation is encoded in a linear system instead of using a linear multi-step method.
Cao et al. \cite{CaoPoisson-2013} present a Hamiltonian simulation algorithm and a scalable quantum circuit design to solve the Poisson equation in $d$ dimensions, where a detailed implementation of the HHL algorithm is also discussed, and the number of quantum operations and the number of qubits used by the circuit are almost linear in $d$ and polylog in $\delta^{-1}$ to produce a quantum state encoding the solution.  Montanaro and Pallister \cite{qFEM-2016} use the quantum linear systems algorithm to solve a linear system resulting from the finite element discretization for the Poisson equation. The state output is then post-processed to compute
a linear functional of the solution. Their algorithms have exponential improvement with respect to $d$ and achieve at most a polynomial speed-up for fixed $d$ due to lower bounds on the cost of post-processing the state.
We remark that a comprehensive review of the literature on quantum differential equations solvers can be found in \cite{Childs-2021}, for example, the quantum algorithms for the wave equation and the hyperbolic equations in \cite{WaveCosta-2019,Fillion-Lorin-2019}.
Other work related to our paper will be reviewed in subsequent sections.

As mentioned above the linear ODEs
\[\frac{\d x}{\d t} = Ax + b, \qquad \mbox{$A$ and $b$ are time-independent}\]
can be simulated by the quantum algorithm in \cite{Berry-Childs-Ostrander-2017} with an exponential improvement in complexity over the quantum difference method in \cite{Berry-2014}. Because of this, it seems that one can first discretize the space variable to obtain a system of ordinary differential equations and then apply the algorithm in \cite{Berry-Childs-Ostrander-2017} to solve the resulting linear system. However, the algorithm requires that $A$ and $b$ are independent of time \cite[Sect.~8]{Berry-Childs-Ostrander-2017}, which may not be applicable to the problems in this paper.
For this reason, in this article we still explore the quantum difference method proposed by Berry in \cite{Berry-2014} to solve several time-dependent problems, analyze in detail the time complexity of the algorithms and make comparisons with the classical treatments.
It should be pointed out that quantum algorithms for differential equations with time-dependent coefficients $A(t)$ and $b(t)$ have been discussed in \cite{Childs-Liu-2020}, where a global approximation based on the
spectral method is employed as an alternative to the more straightforward finite difference method.

Quantum algorithms for linear ODEs and PDEs have been extensively studied, so the focus of this paper is not to develop new and novel quantum algorithms, rather we first try to understand whether  different classical time discretizations, including explicit and implicit schemes, make any difference for quantum algorithms.  The main part of this paper is on a prototype multiscale hyperbolic system with the aim to understand whether one needs to use state-of-the-art multiscale methods for quantum algorithms.

In Section \ref{sec:heat} (resp. Section \ref{sec:1ohyperblic}), we consider the computation of the high-dimensional linear heat equation (resp. linear hyperbolic equation). We find that, for the heat equation, explicit (forward Euler) and implicit (Crank-Nicolson) methods in time have comparable time-complexity for quantum algorithms.  Section \ref{sec:multiscale} discusses the quantum difference method for solving a prototype multiscale problem-a hyperbolic system with stiff relaxation, including three AP schemes and one explicit scheme. The time complexity of both the classical and quantum treatments in the explicit scheme is proportional to $1/\varepsilon$, while the run time of the AP schemes is found to be {\it independent} of the scaling parameter $\varepsilon$, which demonstrates that there is still a need to develop AP schemes for multiscale problems in quantum computing.

One important remark is that most of the quantum algorithms for ODEs and PDEs do {\it not} actually provide the solutions of the equations in classical form. Rather they only prepare the quantum states whose amplitudes encode those solutions. Thus these are sometimes termed subroutines rather than full quantum algorithms. An end-to-end comparison between classical and quantum algorithms would only be possible if one is able to address costs in the preparation of the initial quantum state as well as the final measurements protocol where the classical description of the relevant solutions can be extracted \cite{aaronson2015read}.
While quantum advantages are significant for high dimensional problems, for all problems discussed in low dimensions,  we find that the quantum speedup is at most polynomially scaled with $N_x$ and could be reduced or even evaporated when considering the post-processing measurement steps.

\section{The quantum difference method}\label{sec:QDM}

In this section, we describe the linear system approach introduced in \cite{Berry-2014} for solving linear ordinary differential equations. In view of the classical difference discretization, we term this method as a quantum difference method. For ease of presentation, we take a simple ODE for example. The ODE is
\begin{equation}\label{odef}
\begin{cases}
\frac{{\rm d} u(t)}{{\rm d}t} = - a u(t), \quad t \in (0,1), \\
u(0) = u_0,
\end{cases}
\end{equation}
where $a>0$ is a constant and $u(t)$ is a scalar function of $t$.

\subsection{Quantum difference schemes}

The quantum difference method is proposed in \cite{Berry-2014} for solving the ODEs. Instead of using a linear multi-step method, a truncation of the Taylor series of the propagator for the differential equation is considered in \cite{Berry-Childs-Ostrander-2017}. The same idea in \cite{Berry-2014} is also considered in \cite{Lin-2022} with the forward Euler discretization applied. We are interested in  examining the difference between  explicit and implicit quantum methods. To do so, we use the general $\theta$-scheme.

Let $t_n = n\tau$ for $n = 0,1,\cdots,N_t$, where $\tau$ is the time step and $N_t = 1/\tau$ is the number of time intervals. The numerical solution at $t_n$ is denoted by $u^n$.
The first equation in  \eqref{odef} can be discretized by the $\theta$-scheme as
\[\frac{u^{n+1}-u^n}{\tau} = - a ( \theta u^{n+1} +  (1-\theta) u^n) , \quad n = 0,1,\cdots, N_t-1,\]
where $\theta \in [0,1]$. On a classical computer, one just march in time to   obtain $u^1, \cdots, u^{N_t}$.
The above equation can be rewritten as
\[-(1 - (1-\theta) a \tau) u^n + (1 + \theta a\tau)u^{n+1}   = 0, \quad n = 0,1,\cdots, N_t-1.\]
Let
\[s_1 = 1 - (1-\theta) a \tau, \quad s_2 = 1 + \theta a \tau.\]
The above system is then written in matrix form as
\begin{equation}\label{ODEforwardAxb}
\bb{A} \bb{u} = \bb{b},
\end{equation}
where $\bb{u} = [u^1, \cdots, u^{N_t}]^T$, $\bb{b} = [s_1u^0, 0, \cdots, 0]^T$ and
\[\bb{A} =
\begin{bmatrix}
s_2       &           &            &           &             &       \\
-s_1      &   s_2     &            &           &             &       \\
          &  -s_1     &  s_2       &           &             &       \\
          &           &   \ddots   &   \ddots  &             &       \\
          &           &            &           &     -s_1    &   s_2    \\
\end{bmatrix}_{N_t \times N_t}\,.
\]

The solution of the differential equation at all times will be encoded by using just one state $\ket{\bb{u}}$ corresponding to the solution of \eqref{ODEforwardAxb}. To obtain the quantum speedup, one needs to choose a suitable quantum algorithm to solve the above linear system.

\subsection{Quantum linear systems algorithms}

\subsubsection{Input model}
We first state the quantum linear systems problem.
\begin{definition}[Quantum linear systems problem (QLSP)]
Consider the system of linear equations $Ax = b$, where $A$ is an $N \times N$ Hermitian matrix, and vectors $x = [x_1,\cdots,x_N]^T$ and $b = [b_1,\cdots,b_N]^T$.  Assume that the vectors $x$ and $b$ are encoded as
\[\ket{x} = \frac{1}{N_x} \sum\limits_{i=1}^N x_i \ket{i} \qquad \mbox{and} \qquad
\ket{b} = \frac{1}{N_b} \sum\limits_{i=1}^N b_i \ket{i},\]
where $N_x = (x_1^2+\cdots+x_N^2)^{1/2}$ and $N_b = (b_1^2+\cdots+b_N^2)^{1/2}$ are normalization constants.
The aim of any algorithm to solve QLSP (such an algorithm is called QLSA) is the following. When given access to $A$ and $b$, one aims to prepare a quantum state $\ket{x'}$ that is $\eta$-close to $\ket{x}$, i.e., $\| \ket{x'} - \ket{x} \|  \le \eta$.
\end{definition}

Let's state the concept of time complexity.
On a classical computer, the time complexity of an algorithm is usually a count of the number of basic operations of addition and multiplication. For example, when comparing the time complexity of their proposed quantum linear solver in \cite{HHL-2009} with that of the classical matrix inversion algorithms, the authors state that the conjugate gradient (CG) method uses $\mathcal{O}(\sqrt{\kappa} \log (1/\delta) )$ matrix-vector multiplications each taking time $\mathcal{O}(Ns)$ for a total running time of $\mathcal{O}(Ns\sqrt{\kappa} \log (1/\delta) )$ for a positive definite matrix, where $s$ is the sparsity number, meaning it has at most $s$ nonzero entries per row, $\kappa$ is the condition number of the matrix, and $\delta$ is the expected error bound of the algorithm.

On a quantum computer, the time complexity of an algorithm can be measured by the gate complexity which refers to the total number of 2-qubit gates used in the algorithm \cite{Childs-2017}. This is because quantum algorithms are usually represented by quantum circuits, whose basic operation is a 2-bit quantum logic gate.
The more commonly used measure of the time complexity is the query complexity \cite[Sect.~1.10]{Lin-2022}.
For quantum linear systems problems, the query complexity is usually measured in terms of the number of calls an oracle, for instance to entries of a sparse matrix $A$ as in \cite{Childs-2017}. The sparse matrices can also be block encoded using standard methods \cite{Low-2019}. In this case one can consider the number of queries to the oracles involved in the block encoding of the coefficient matrix $A$ \cite{Costa-An-Sanders-2021,Lin-2022}. We assume $\| A \|_{\text{max}} < 1$ throughout the discussion, unless otherwise stated, since the complexity can have a contribution proportional to $\| A \|_{\text{max}}$ \cite{Jin-Liu-2022}. Otherwise we can simply replace it by the re-scaled matrix $A/\alpha$ for some $\alpha > \| A \|_{\text{max}}$, where $\| A \|_{\text{max}} = \max_{ij} |A_{ij}|$.

In this article we are concerned with the sparse access to the matrix since the resulting coefficient matrices for PDEs are usually sparse, defined in the following way \cite{Berry-Childs-Ostrander-2017,Jin-Liu-2022}.

\begin{definition}[Input model for the matrix]
Let $A$ be a Hermitian matrix with the $(i,j)^{\text{th}}$ entry denoted by $A_{ij}$. Sparse access to $A$ is referred to as a 4-tuple $(s, \| A \|_{\text{max}}, O_A, O_{F})$. Here, $s$ is the sparsity of $A$; $O_A$ is a unitary black box which can access the matrix elements $A_{ij}$ such that
\[O_A |j\rangle|k\rangle|z\rangle = |j\rangle|k\rangle|z\oplus A_{jk}\rangle\]
for any $j,k \in \{1,2,\cdots,N\}=:[N]$, where the third register holds a bit string representing of $A_{jk}$;
$O_F$ is a unitary black box which allows to perform the map
\[O_F |j\rangle|l\rangle = |j\rangle|F(j,l)\rangle\]
for any $j\in [N]$ and $l \in [s]$, where the function $F$ outputs the column index of the $l^{\text{th}}$ non-zero elements in row $j$.
\end{definition}

\begin{definition}[Input model for the initial state]
Let $b$ be a vector. Access to $b$ is referred to as a unitary $U_{initial}$  such that $U_{initial}\ket{0}=\ket{b}$. See \cite{gleinig2021efficient, zhang2022quantum} for example.
\end{definition}

With these definitions, the query complexity denotes the number of times oracles $O_A$, $O_F$ and $U_{initial}$ are used throughout the protocol.
As reviewed later, we shall apply the optimal QLSA proposed in \cite{Costa-An-Sanders-2021} to solve the linear system \eqref{ODEforwardAxb}, in which the query complexity is in terms of calls to a block encoding of the coefficient matrix, rather than the more fundamental oracles for positions of nonzero entries of sparse matrices, defined below \cite{Jin-Liu-2022,Lin-2022}.

\begin{definition}[Block access to the matrix]
Let $A$ be a $m$-qubit Hermitian matrix, $\delta_A>0$ and $n_A$ is a positive integer. An $(m+n_A)$-qubit unitary matrix $U_A$ is a $(\alpha_A, n_A, \delta_A)$-block encoding of $A$ if
\[\|A-\alpha_A\langle 0^{n_A}|U_A |0^{n_A} \rangle \|\le \delta_A.\]
Block access to $A$ is then the 4-tuple $(\alpha_A, n_A, \delta_A, U_A)$ where $U_A$ is the unitary black-box block-encoding of $A$.
\end{definition}

It is possible to create block access to $A$ from sparse access to $A$ \cite{Jin-Liu-2022,Lin-2022,Low-2019}. If standard methods \cite{Low-2019} are used to construct the block access from the sparse access, there will be a multiplicative factor $s$ in the query complexity of block encoding \cite{Costa-An-Sanders-2021}.

The coefficient matrix $A$ in a linear system is prepared on a classical computer and stored in an external database with the bit string representing of $A_{jk}$. Then one can construct the matrix in a quantum computer by using the input model to query the non-zero entries through the external database.

\subsubsection{Review of quantum linear systems algorithms}

Ref. \cite{Berry-2014} uses the well-known HHL algorithm for solving \eqref{ODEforwardAxb}. This algorithm is the earliest quantum algorithm for solving systems of linear equations, proposed by Harrow, Hassidim and Lloyd in 2009 in \cite{HHL-2009}. The gate complexity of this algorithm is $\mathcal{O}(\log(N) s^2\kappa^2 /\delta )$.
 When the condition number $\kappa$ of the matrix is not very large, the algorithm has an exponential speedup over the classical algorithm. In fact, if the matrix is positive definite, then the time complexity of the CG method is
$\mathcal{O}(N s\sqrt{\kappa} \log (1/\delta) $, which exhibits exponential growth in the order of the matrix compared to the HHL algorithm.

However, the direct use of the HHL algorithm does not necessarily yield quantum acceleration advantages when solving the linear system of equations obtained by numerically discretizing a partial differential equation. On the one hand, the condition number of the coefficient matrix of the induced linear system increases with mesh refinement; on the other hand, the time complexity of the HHL algorithm depends on $1/\delta$, and $\delta$ should be smaller than or in the same order as the error induced by the numerical discretizations. These two observations may offset the quantum advantage in terms of matrix order. For example, for the five-point difference scheme of the two-dimensional Poisson equation \cite{CaoPoisson-2013,Lin-2022}, the condition number of its coefficient matrix is $\kappa = \mathcal{O}(h^{-2})$, and the accuracy in maximum norm is $\mathcal{O}(h^2)$, where $h=1/N_x$ is the mesh size in the direction of $x$ or $y$. Let $\delta = h^2$. One has
\[C_{\text{CG}} = \mathcal{O}(\log(N_x) N_x^2), \qquad
Q_{\text{HHL}} = \mathcal{O}(\log(N_x) N_x^6 ),\]
since $N = \mathcal{O}(N_x^2)$ and $s=\mathcal{O}(1)$. Here and later, the symbols $C$ and $Q$ are used to denote the time complexity for the classical and quantum algorithms, respectively.
As observed, the time complexity of the classical method is {\it smaller} than that of the quantum method. In fact, for the PDE problem in lower dimensions, the more obvious contribution to the increase in complexity of the HHL algorithm is the term associated with $\delta$, because it increases exponentially compared with the CG method. Of course, as long as the order of the matrix increases faster than the condition number and $1/\delta$, the quantum speedup can be obtained. This can be seen in higher dimensional problems, such as the high-dimensional heat equation discussed in this paper.

To improve the accuracy dependence, Childs et al. \cite{Childs-2017} replaces the Hamiltonian simulation with Fourier method or Chebyshev method to obtain $\text{polylog} (1/ \delta)$ dependence, which is similar to that of the classical method, both of which have exponential acceleration in accuracy over the HHL algorithm.  Childs also uses their quantum method to solve some linear partial differential equations in \cite{Childs-2021}. It is worth pointing out that Childs et al. also exploited the variable-time amplitude amplification (VTAA) proposed by Ambainis in \cite{Ambainis-2012} to improve the conditional number dependence in \cite{Childs-2017}, reducing the $\kappa$-dependence from quadratic to nearly linear with the query complexity given by $\mathcal{O}(\kappa \text{polylog}(\kappa/\delta))$.
However, the VTAA procedure is highly complicated, making it challenging to implement in practice due to the multiple rounds of recursive amplitude amplifications \cite{Costa-An-Sanders-2021}, and it is still asymptotically sub-optimal by a factor of $\log(\kappa)$. To address this issue, alternative approaches based on adiabatic quantum computing (AQC) have been developed in recent years. One can refer to \cite{Costa-An-Sanders-2021} for a comprehensive review of the literature along this line, where the optimal scaling with the condition number is achieved with query complexity $\mathcal{O}(\kappa \log (1/\delta))$ by using a discrete quantum adiabatic theorem proved in \cite{Dranov-Kellendonk-Seiler-1998}, which completely avoids the heavy mechanisms of VTAA or the truncated Dyson-series subroutine from previous methods related to the AQC \cite{Subasi-Somma-2019}. Unlike previous algorithms like HHL and CG, these algorithms assume access to block-encodings of the matrix $A$ instead, thus do not carry dependency the sparsity of $A$.

It is known that a quantum algorithm must make at least $\Omega(\kappa \log(1/\delta))$ queries in general to solve the sparse quantum linear system problems, where the notation $f = \Omega(g)$ means $g = \mathcal{O}(f)$ \cite{An-Lin-2022,Costa-An-Sanders-2021}. Therefore, the method in \cite{Costa-An-Sanders-2021} is already optimal in the scaling with precision  $\delta$ and  condition number $\kappa$. Another method to obtain the optimal complexity is the quantum singular value transformation (QSVT) \cite{Gilyen-Su-Low-2019}, which can solve quantum linear system problems for general matrices without the need of dilating the matrix into a Hermitian matrix. The query complexity for solving $Ax=b$ is $\mathcal{O}(\kappa^2/\xi \log (\kappa/(\xi \delta))$, where $\xi = \|A^{-1}b\|$. In the best case scenario that $b$ has an $\Omega(1)$ overlap with the left-singular vector of $A$ with respect to the smallest singular value,  $\xi = \Omega(\kappa)$ and hence the optimal run time $\mathcal{O}(\kappa \log (1/\delta))$ is obtained \cite{Lin-2022}.

Instead of using the HHL algorithm as in \cite{Berry-2014}, we intend to apply the optimal QLSA proposed in \cite{Costa-An-Sanders-2021} to solve the linear system \eqref{ODEforwardAxb}. Note that the query complexity in \cite{Costa-An-Sanders-2021} is in terms of calls to a block encoding of the coefficient matrix, rather than the more fundamental oracles for positions of nonzero entries of sparse matrices. In this article we are concerned with the sparse access to the matrix since the resulting coefficient matrices for PDEs are usually sparse. In this case, there will be a multiplicative factor of $s$ in the sparsity $s$ if standard methods \cite{Low-2019} are used to construct the block access from the sparse access. Therefore, the query complexity with respect to the sparse access to matrices can be written as
\begin{equation}\label{cpCAS}
Q = \mathcal{O} ( s \kappa \log(1/\delta)   ).
\end{equation}
On the other hand, the gate complexity may be quantified by $\mathcal{O}( Q \text{poly} (\log Q, \log N) )$, which is larger than the query complexity only by logarithmic factors \cite{Childs-2017,Lin-Tong-2020,Costa-An-Sanders-2021}.

Note that both the HHL algorithm and the optimal algorithm in \cite{Costa-An-Sanders-2021} are for Hermitian matrices. If $\bb{A}$ is not Hermitian, we then consider
\[\bb{H} = \begin{bmatrix} \bb{O} & \bb{A} \\  \bb{A}^T  & \bb{O}  \end{bmatrix},\]
where $\bb{A}$ is confined to be real matrix throughout the paper.

\subsection{Complexity analysis}

\subsubsection{Query complexity of solving the resulting QLSP}

Applying the elementary transform of block matrices, one easily gets
\[{\rm det} (\lambda \bb{I} - \bb{H} ) = {\rm det} (\lambda^2 \bb{I} - \bb{A}\bb{A}^T ).\]
Hence the eigenvalues of $\bb{H}$ are $\lambda = \pm \sigma_1, \pm \sigma_2, \cdots, \pm \sigma_{N_t}$, where $ \sigma_i$ is the singular values of $\bb{A}$ for $i = 1,\cdots, N_t$.
We are able to derive the time complexity of the quantum difference method described as follows.

\begin{theorem}
Let $a>0$ be a constant and the time step $\tau$ satisfies $\tau < 1/ (a(1-\theta))$. Then the condition number and the sparsity of $\bb{H}$ satisfy
\[ \kappa =  \mathcal{O}(N_t) \qquad \mbox{and} \qquad  s = \mathcal{O}(1).\]
The time complexities of solving the ODE problem \eqref{odef} by using the classical and quantum difference methods are $C = \mathcal{O}(N_t)$ and $Q = \mathcal{O}( N_t \log(1/\delta))$, respectively.
\end{theorem}
\begin{proof}
  (1)   Applying the Gershgorin circle theorem \cite{Horn2013} to $\bb{A}\bb{A}^T$, one easily finds the eigenvalues of $\bb{H}$ satisfy
  \[|\lambda|_{\min} \ge a\tau, \quad |\lambda|_{\max} \le 2 + |2\theta-1| a\tau.\]
  Since the calculation is relatively simple, we omit the details.

  (2) The time complexity of the classical methods simply involves counting the number of basic operations of addition and multiplication. At each iteration step, the number of operations is $\mathcal{O}(1)$, thus the total time complexity is $C = \mathcal{O}(N_t)$.

  (3) The matrix $\bb{H}$ has sparsity number $s=2$ and condition number
  \[\kappa \le \frac{2 + |2\theta-1| a\tau}{a\tau} \lesssim \frac{1}{\tau} = N_t.\]
  Then $Q$ is obtained by plugging these quantities in \eqref{cpCAS}.
\end{proof}

The error of the Euler method ($\theta=0$) is $\mathcal{O}(\tau)$ in $L^\infty$ norm, and that of the C-N method ($\theta=1/2$) is $\mathcal{O}(\tau^2)$. Let the step size of the C-N method be $\tau = 1/N_t$. Then to obtain the same error bound $\delta = \mathcal{O}(\tau^2)$, the step size of the Euler method has to be taken as $\tau^2 = \mathcal{O}(1/N_t^2)$.
Thus, for the classical methods we have
\[C_{\text{Euler}} = \mathcal{O}(N_t^2), \qquad C_{\text{CN}} = \mathcal{O}(N_t). \]
And the corresponding result for the quantum difference method is
\[Q_{\text{Euler}} = \mathcal{O}(N_t^2\log N_t), \qquad Q_{\text{CN}} = \mathcal{O}(N_t \log N_t). \]

For this simple first-order ODE problem, the quantum difference method does not gain quantum advantages. The complexities are only comparable to the ones of the classical methods if the logarithmic terms are neglected.
As mentioned in the introduction, however, the linear ODEs
\[\frac{\d x}{\d t} = Ax + b, \qquad \mbox{$A$ and $b$ are time-independent}\]
can be simulated by the quantum algorithm in \cite{Berry-Childs-Ostrander-2017} with an exponential improvement in time complexity over the quantum difference method in \cite{Berry-2014}. In this case, the query complexity is only $\mathcal{O}(1)$, hence the quantum speedup is recovered. It should be noted that the exact solution of \eqref{odef} is simply $u(t) = \e^{-a t}$, which also requires only $\mathcal{O}(1)$ run time on a classical computer.

\subsubsection{Post-processing: computation of physical quantities of interest}

 The solution of the quantum linear systems problem actually corresponds to the `history-state' solution of the
differential equations, which is the quantum state that is a superposition of the solution at all temporal and spatial points.
However, measurements on the quantum `history state' are required to read out the classical solutions. When measurement costs are included, the quantum speedup can be in cases greatly reduced or even disappear entirely \cite{aaronson2015read}.

When solving ODEs and PDEs, the actual desired outcomes of the problem are physical quantities of interest that are  associated with the solutions of the ODEs and PDEs. The selection of these physical observables depends on the problem under consideration. For example, Ref.~\cite{Lin-Montanaro-Shao-2020} computes the amount of heat $\int_S u $ in a given region $S$ when solving the heat equation, and Ref.~\cite{qFEM-2016} approximates the inner product $\langle u, r \rangle$ between the solution $u$ and a fixed function $r$. We also refer the reader to \cite{Jin-Liu-2022} for how to compute ensemble averages of physical observables for nonlinear PDEs using quantum algorithms.

Let $O$ be an observable with $\mu:=\langle O \rangle = \langle \psi| O | \psi \rangle$ being the expectation value, where $\ket{\psi}$ is a quantum state. Suppose that we conduct $n$ experiments with the outcomes recorded as $\mu_1, \cdots, \mu_n$. By the law of large numbers, we have
\[{\rm P}_{\text{r}}\Big( \Big|\frac{\mu_1+\cdots + \mu_n}{n} -  \mu\Big| < \varepsilon \Big) \ge 1 - \frac{\text{Var}(O)}{n\varepsilon^2},\]
where $\text{Var}(O)$ is the variance. For a given lower bound $p$, the number of samples required to estimate $\langle O \rangle$ to additive precision $\varepsilon$ satisfies
\[1 - \frac{\text{Var}(O)}{n\varepsilon^2} \ge p \quad \Longrightarrow \quad
n \ge \frac{1}{1-p} \frac{\text{Var}(O)}{\varepsilon^2}.\]
This implies a multiplicative factor $\text{Var}(O)/\varepsilon^2$ in the total time complexity \cite{Lin-2022}. We note that in a block-encoding scheme we can also use amplitude estimation to enhance the error scaling to $O(1/\varepsilon)$ (e.g. \cite{rall2020quantum, Jin-Liu-2022}). In this paper we don't focus on this scaling and use the simpler protocol with $O(1/\varepsilon^2)$ scaling since our goal is only to compare quantum algorithms with different discretisation schemes.

We remark that the variance $\text{Var}(O)$ may scale as $\mathcal{O}(N_t^2) = \mathcal{O}(\tau^{-2})$ for the QLSA. To see this, let us denote the quantum state of the solution to the QLSA by
\[\ket{\widetilde{\bb{u}}} = [\widetilde{\bb{u}}^1; \cdots; \widetilde{\bb{u}}^{N_t}], \qquad
\widetilde{\bb{u}}^n = \frac{1}{N_u}\bb{u}^n,\]
where $\bb{u}^n$ is the solution vector at time $t=t_n$ and the normalization constant is
\[N_u = \|\bb{u}\| = (\|\bb{u}^1\|^2+ \cdots + \|\bb{u}^{N_t}\|^2)^{1/2}.\]
For simplicity, we assume that $\|\bb{u}^1\| = \cdots = \|\bb{u}^{N_t}\| = \|\bb{u}^0\|$ and hence $N_u = \sqrt{N_t} \|\bb{u}^0\|$. Let $O_{\bb{i}} = \ket{\bb{i}}\bra{\bb{i}}$, $O_n = \ket{n}\bra{n}$ and
\[O_{\bb{i}}^n = O_n \otimes O_{\bb{i}} = \ket{n,\bb{i}}\bra{n,\bb{i}},\]
where $\ket{n}$ is of size $N_t$. We consider a simple observable
\[\rho(t = t_n, x_{\bb{i}}) := (\bb{u}^n)^\dag  O_{\bb{i}} \bb{u}^n
= \bb{u}^\dag (O_n \otimes O_{\bb{i}}) \bb{u}
= N_t N_{u_0}^2 \cdot \langle \widetilde{\bb{u}} | O_{\bb{i}}^n | \widetilde{\bb{u}} \rangle.\]
The expectation $\langle O_{\bb{i}}^n \rangle := \langle \widetilde{\bb{u}} | O_{\bb{i}}^n | \widetilde{\bb{u}} \rangle$ satisfies the condition that ${\rm Var}( O_{\bb{i}}^n)$ is bounded.
In this case, however, we must evaluate $\langle O_{\bb{i}}^n \rangle$ to precision $\mathcal{O}(\varepsilon/(N_tN_{u_0}^2))$, which increases the number of samples by another factor $(N_tN_{u_0}^2)^2$. Note that we cannot resolve the issue by considering $N_t N_{u_0}^2 O_{\bb{i}}^n$ directly since ${\rm Var}( N_tN_{u_0}^2 O_{\bb{i}}^n ) = (N_tN_{u_0}^2)^2 {\rm Var}( O_{\bb{i}}^n )$ gives the same factor.

 A simple way to overcome this problem has been addressed in the original paper \cite{Berry-2014} by adding $N_t$ copies of the final state $\bb{u}^{N_t}$. That is, we add the following additional equations
\begin{equation}\label{dilation}
\bb{u}^{n+1} - \bb{u}^n = 0, \quad n = N_t, \cdots, 2N_t,
\end{equation}
which is referred to as the dilation procedure \cite{Lin-2022}. In fact, let the padded state vector be
\begin{equation} \label{padding}
\widehat{\bb{u}}
 = [\widehat{\bb{u}}^1; \cdots; \widehat{\bb{u}}^{N_t};\widehat{\bb{u}}^{N_t},\cdots, \widehat{\bb{u}}^{N_t}  ]
  = \ket{0} \otimes \bb{x} + \ket{1} \otimes \bb{y},
\end{equation}
where $\ket{0} = [1,0]^T$, $\ket{1} = [0,1]^T$, and the unnormalized vectors are
\[\bb{x} = [\widehat{\bb{u}}^1; \cdots; \widehat{\bb{u}}^{N_t}], \qquad
\bb{y} = [\widehat{\bb{u}}^{N_t},\cdots, \widehat{\bb{u}}^{N_t}  ].\]
Noting that
\[\|\widehat{\bb{u}}^n\|^2 = \frac{1}{2N_t} = \frac{1}{2N_t N_{u_0}^2} \|\bb{u}^n\|^2, \qquad n = 1,\cdots, N_t,\]
we further define
\[\widehat{O} = {\rm diag}(O_{\bb{i}}, \cdots, O_{\bb{i}}) = I_{N_t} \otimes O_{\bb{i}},\]
where $I_{N_t}$ is the identity matrix with order $N_t$, and obtain
\begin{align*}\rho(t = t_{N_t} ,x_{\bb{i}})
& = (\bb{u}^{N_t})^\dag  O_{\bb{i}}  \bb{u}^{N_t}
= 2 N_tN_{u_0}^2 \cdot (\widehat{\bb{u}}^{N_t})^\dag  O_{\bb{i}}  \widehat{\bb{u}}^{N_t} \\
& = 2 N_{u_0}^2\bb{y}^\dag  \widehat{O}  \bb{y} = 2 N_{u_0}^2 \langle \widehat{\bb{u}} | O | \widehat{\bb{u}} \rangle,
\end{align*}
where
\[O = {\rm diag}( \bb{O}, \cdots, \bb{O}, O_{\bb{i}}, \cdots, O_{\bb{i}}), \quad \bb{O}~~\mbox{is the zero matrix}. \]
Obviously, ${\rm Var}( O )$ is bounded. It's worth pointing out that the solution vector in \eqref{padding} only requires one ancilla qubit, and one can directly evaluate $\bb{y}^\dag  \widehat{O}  \bb{y}$ by measuring the ancilla qubit and obtain 1.

We remark that the quantum solution, proportional to $\bb{u}$, is a history state. In order to recover the solution at a time $N_t$, one needs to project onto the solution at $N_t$. However, since the solution of Eq.~\eqref{odef} decays exponentially in time, the success probability of such a projection is exponentially small. This problem exists for all QLSAs that outputs a history state.

However, this probability can be raised if we apply amplitude amplification to this scheme. In this case we can raise the probability to $\Omega(1)$ with $\mathcal{O}(g)$ repetitions of the QLSA, where $g = \max_{t \in [0,T]} \| \bb{u}(t) \| /  \| \bb{u}(T) \|$ characterizes the decay of the final state relative to the initial state. This implies a new factor $g$ in the final time complexity.
We remark that the parameter $g$ has been included in the complexities of the quantum algorithms in \cite{Berry-Childs-Ostrander-2017,Childs-Liu-2020}. For convenience, in the subsequent discussion, we ignore this parameter in the time complexity.

\section{The linear heat equation}\label{sec:heat}

Consider the following initial-boundary value problem of linear heat equation
\begin{equation}\label{eq:heat}
\begin{cases}
& u_t(\bb{x},t) = \Delta u(\bb{x},t) \quad  \text{in}~~ \Omega :=(0,1)^d, \quad 0<t<1, \\
& u(\bb{x},0) = u_0(\bb{x}), \\
& u(\cdot, t) = 0 \quad \text{on}~~\partial \Omega,
\end{cases}
\end{equation}
where $u_0(\bb{x})$ is sufficiently smooth.

During the revision of our paper, we found that the heat equation was also considered in \cite{Lin-2022} with a forward Euler discretization in time applied. In addition, the authors in \cite{Lin-Montanaro-Shao-2020} studied in detail the complexities of ten classical and quantum algorithms for solving the heat equation in the sense of approximately computing the amount of heat in a given region, in which the quantum linear equations method was discussed in Theorem 17 there.
However, from the numerical analysis point of view it is of interest to know if  some discretization ways are preferred to others for the quantum difference methods. Our study reveals that  different discretizations  (forward Euler and Crank-Nicolson) don't make a difference in the quantum scenario. This is a phenomenon not reported  before to our knowledge.

\subsection{The finite difference schemes}

For simplicity we only provide the detail for the one-dimensional case. Consider $N_t+1$ steps in time $0 = t_0<t_1<\cdots<t_{N_t} = 1$ and $N_x+1$ spatial mesh points $0<x_0<x_1<\cdots<x_{N_x} = 1$ by setting $t_n = n\tau$ and $x_j = jh$, where $\tau = 1/N_t$ and $h = 1/N_x$. The central difference discretization gives
\[\frac{{\rm d}}{{\rm d}t}u_j(t) = \frac{u_{j-1}(t) - 2u_j(t) + u_{j+1}(t)}{h^2},  \quad j = 1,\cdots, N_x-1.  \]
Let $\bb{u}(t) = [u_1(t), \cdots, u_{N_x-1}(t)]^T$. One has
\begin{equation}\label{odeheat}
\frac{{\rm d}}{{\rm d}t} \bb{u}(t)  =  \frac{1}{h^2} \bb{L}_h \bb{u}(t) + \frac{1}{h^2} \bb{b}(t),
\end{equation}
where,
\[\bb{L}_h =
\begin{bmatrix}
-2  &  1       &           &      &    \\
 1  & -2       & \ddots    &      &    \\
    &  \ddots  & \ddots    &  \ddots    &    \\
    &          & \ddots    & -2   & 1  \\
    &          &           &  1   & -2 \\
\end{bmatrix}_{(N_x-1) \times (N_x-1)}
, \qquad
\bb{b}(t) =
\begin{bmatrix}
u_0(t) \\
0  \\
\vdots\\
0 \\
u_{N_x}(t) \\
\end{bmatrix}.
\]
The time is discretized by using the $\theta$-scheme as
\[\frac{\bb{u}^{n+1} - \bb{u}^n}{\tau} = \theta \frac{1}{h^2} \bb{L}_h \bb{u}^{n+1} + (1-\theta) \frac{1}{h^2} \bb{L}_h \bb{u}^n+ \theta \frac{1}{h^2} \bb{b}^{n+1} + (1-\theta) \frac{1}{h^2} \bb{b}^n,\]
which can be marched forward in time for $\theta=0$ or by solving a linear system for $\theta\in (0, 1]$  on a classical computer.

Let $\beta = \tau/h^2$. The above equation can be written as
\begin{equation}\label{heatC}
-\bb{B}\bb{u}^n + \bb{A} \bb{u}^{n+1} = \bb{f}^{n+1},
\end{equation}
where
\[\bb{A} = \bb{I}- \theta \beta \bb{L}_h, \quad  \bb{B} = \bb{I} + (1-\theta) \beta \bb{L}_h, \quad
\bb{f}^{n} = \theta\beta\bb{b}^{n+1} + (1-\theta) \beta \bb{b}^n.\]
By introducing the notation $\bb{U} = [\bb{u}^1; \cdots ;\bb{u}^{N_t}]$, where ``;" indicates the straightening of $\{\bb{u}^i\}_{i\ge 1}$ into a column vector, and one  obtains the following linear system
\begin{equation}\label{heatLF}
\bb{L} \bb{U} = \bb{F},
\end{equation}
where
\[\bb{L} =
\begin{bmatrix}
\bb{A}  &            &           &            \\
-\bb{B} & \bb{A}     &           &            \\
        &\ddots      & \ddots    &    \\
        &            & -\bb{B}   & \bb{A}     \\
\end{bmatrix}
, \qquad
\bb{F} =
\begin{bmatrix}
\bb{f}^1 + \bb{B}\bb{u}^0 \\
\bb{f}^2  \\
\vdots\\
\bb{f}^{N_t} \\
\end{bmatrix}.
\]
Note that the hermitian matrix for the quantum linear system problem is
\[\bb{H} = \begin{bmatrix} \bb{O} & \bb{L} \\  \bb{L}^T  & \bb{O}  \end{bmatrix}.\]

The matrix $\bb{L}_h$ of the $d$-dimensional problem will be replaced by
\[\bb{L}_{h,d} = \underbrace{\bb{L}_h\otimes \bb{I} \otimes \cdots \otimes \bb{I}}_{d~\text{matrices}} + \bb{I} \otimes \bb{L}_h\otimes \cdots \otimes \bb{I} + \cdots + \bb{I} \otimes \bb{I} \otimes \cdots \otimes \bb{L}_h,\]
and everything else remains the same. In addition, we only consider solving the high-dimensional heat equation with homogeneous Dirichlet boundary condition. Other boundary conditions can be similarly treated with possible different condition number $\kappa$ of the coefficient matrix.

\subsection{Comparison of the time complexity}

For the explicit scheme, i.e. $\theta=0$, the stability condition is $d \tau/h^2 \le 1$ and the error in $L^\infty$ norm is $\mathcal{O}(\tau + dh^2)$. We thus choose $\tau \le h^2/(4d)$ for $\theta=0$. The C-N scheme is unconditionally stable, and its error is $\mathcal{O}(\tau^2 + dh^2)$.
In view of the constraint in the following theorem, we set $\tau = h/(8d)$ for $\theta=1/2$.

In the following we take the same mesh size $h = 1/N_x$, and set $\delta = dh^2 = d/N_x^2$.
Then,
\begin{itemize}
  \item Forward ($\theta=0$): $\tau = \mathcal{O}(h^2/d)$ or $N_t = \mathcal{O}(d N_x^2)$, and $\beta = \mathcal{O}(1/d)$;
  \item C-N ($\theta=1/2$): $\tau = \mathcal{O}(h/d)$ or $N_t = \mathcal{O}(d N_x)$, and $\beta = \mathcal{O}(N_x/d)$.
\end{itemize}

\begin{theorem}
Let $\theta=0$ or $1/2$ and $\tau \le 1/(8d)$. For the explicit scheme, i.e., $\theta=0$, we further assume that the parabolic CFL condition $\beta = \tau/h^2 \le 1/(4d)$ holds. With the above settings, one has
\begin{enumerate}[(1)]
  \item The condition number and the sparsity of $\bb{H}$ satisfy
\[\kappa = \mathcal{O}(N_x^2) \qquad \mbox{and} \qquad s = \mathcal{O}(d).  \]
  \item The time complexity of the classical difference methods for solving the $d$-dimensional heat equation is
  \[C = \begin{cases}
  \mathcal{O}(d^2 N_x^{d+2}) , \quad & \theta = 0, \\
  \mathcal{O}(d^2 N_x^{d+1.5}), \quad & \theta = 1/2.
  \end{cases}\]
  \item The time complexity for the quantum difference methods is
  $Q = \mathcal{O} ( d N_x^2 \log(N_x^2/d)  )$.
\end{enumerate}
\end{theorem}
\begin{proof}

(1) Let $\lambda$ be the eigenvalue of $\bb{H}$, i.e.,
\[\text{det}(\lambda\bb{I}-\bb{H}) = \text{det}(\lambda^2\bb{I}-\bb{L}\bb{L}^T) =:\text{det}(\mu\bb{I}-\bb{L}\bb{L}^T),\]
where $\mu = \lambda^2$ is the eigenvalue of $\bb{L}\bb{L}^T$, or $\mu^{1/2}$ is the singular value of $\bb{L}$. Obviously $\bb{A}$ and $\bb{B}$ are symmetric matrices, and the direct calculation gives
\[\bb{L}\bb{L}^T =
\begin{bmatrix}
\bb{A}^2  &    -\bb{AB}             &                   &            \\
-\bb{BA}  & \bb{B}^2 + \bb{A}^2     & \ddots            &            \\
          &\ddots                   & \ddots            &   -\bb{AB} \\
          &                         & -\bb{BA}          & \bb{B}^2 + \bb{A}^2    \\
\end{bmatrix}.
\]

We again establish the upper and lower bounds of the eigenvalues of $\bb{L}\bb{L}^T$ by using the Gershgorin circle theorem. Consider the 1-D case. Let $\bb{P}^{-1}\bb{L}_h\bb{P} = \Lambda$, where $\Lambda$ is the diagonal matrix consisting of the eigenvalues of $\bb{L}_h$. For convenience, we denote the relation as $\bb{L}_h \sim \Lambda$ and obtain
\[
 \bb{A} \sim \bb{I} - \theta\beta\Lambda = : \Lambda_A,   \qquad  \bb{B} \sim \bb{I} + (1-\theta)\beta\Lambda = : \Lambda_B, \qquad  \bb{AB} \sim \Lambda_A \Lambda_B.
\]
Let $\widetilde{\bb{P}} = \text{diag}(\bb{P}, \cdots, \bb{P})$. Then
\[
\widetilde{\bb{P}}^{-1}(\bb{L}\bb{L}^T)\widetilde{\bb{P}} =
\begin{bmatrix}
\Lambda_A^2  &    -\Lambda_A \Lambda_B             &                   &            \\
-\Lambda_A \Lambda_B  & \Lambda_A^2 + \Lambda_B^2     & \ddots            &            \\
          &\ddots                   & \ddots            &   -\Lambda_A \Lambda_B \\
          &                         & -\Lambda_A \Lambda_B          & \Lambda_A^2 + \Lambda_B^2    \\
\end{bmatrix}.
\]
Noting that the similarity transformation does not change the eigenvalues, thus one can apply the Gershgorin circle theorem to estimate the eigenvalues for this matrix, which implies that the minimum and maximum eigenvalues satisfy
\begin{align*}
&\mu_{\min} \ge \min_{i} \{ \lambda_{A,i}^2 - |\lambda_{A,i}\lambda_{B,i}|, \quad \lambda_{A,i}^2 + \lambda_{B,i}^2 - 2|\lambda_{A,i}\lambda_{B,i}| \},\\
&\mu_{\max} \le \max_{i} \{ \lambda_{A,i}^2 + \lambda_{B,i}^2 + 2|\lambda_{A,i}\lambda_{B,i}| \}.
\end{align*}

The eigenvalues of $\bb{L}_h$ are
\begin{equation}\label{eigLh}
\lambda_{h,i} = -4\sin^2\frac{i\pi}{2N_x} = -4\sin^2\frac{i\pi h}{2}, \quad i = 1,\cdots, N_x-1,
\end{equation}
satisfying $8h^2 \le |\lambda_{h,i}| \le 4$. Since $\lambda_{h,i} < 0$,
\[\lambda_{A,i} = 1 - \theta \beta \lambda_{h,i} > 0.\]
If $\lambda_{B,i} \ge 0$, then
\begin{equation}\label{ABmimus}
|\lambda_{A,i}| - |\lambda_{B,i}|
= 1 - \theta \beta \lambda_{h,i} - 1 - (1-\theta) \beta \lambda_{h,i} = -\beta \lambda_{h,i} \ge 8 \tau >0
\end{equation}
and
\begin{align*}
|\lambda_{A,i}| + |\lambda_{B,i}|
& = 1 - \theta \beta \lambda_{h,i} + 1 + (1-\theta) \beta \lambda_{h,i} \\
& = 2 + (1-2\theta) \beta \lambda_{h,i} =\begin{cases}
2 + \beta \lambda_{h,i} \le 3, \quad & \theta = 0, \\
2, \quad & \theta = 1/2.
\end{cases}
\end{align*}
If $\lambda_{B,i} < 0$, then
\begin{align}
|\lambda_{A,i}| - |\lambda_{B,i}|
& = 1 - \theta \beta \lambda_{h,i} + 1 + (1-\theta) \beta \lambda_{h,i} \nonumber \\
& = 2 + (1-2\theta)\beta \lambda_{h,i} = \begin{cases}
2 + \beta \lambda_{h,i} \ge 1, \quad & \theta = 0, \\
2, \quad & \theta = 1/2,
\end{cases} \label{ABmimus1}
\end{align}
where the condition $\beta\le 1/4$ is used for $\theta=0$, and
\begin{align*}
|\lambda_{A,i}| + |\lambda_{B,i}|
& = 1 - \theta \beta \lambda_{h,i} - 1 - (1-\theta) \beta \lambda_{h,i} \\
& = -\beta \lambda_{h,i} \le \begin{cases}
1, \quad & \theta = 0, \\
4 \beta, \quad & \theta = 1/2.
\end{cases}
\end{align*}

We first bound the maximum eigenvalue. According to the previous calculations, one has
\[\lambda_{A,i}^2 + \lambda_{B,i}^2 + 2|\lambda_{A,i}\lambda_{B,i}|
= ( |\lambda_{A,i}| + |\lambda_{B,i}|)^2
\le \begin{cases}
9, \quad & \theta = 0,\\
\max\{4, (4\beta)^2\}, & \theta = 1/2,\\
\end{cases} \]
which implies
\[\mu_{\max} \le \begin{cases}
9, \quad & \theta = 0, \\
\max\{4, (4\beta)^2\}, & \theta = 1/2.
\end{cases}\]
For the minimum eigenvalue, we obtain from \eqref{ABmimus} and \eqref{ABmimus1} that
\[(\lambda_{A,i}^2 + \lambda_{B,i}^2 - 2|\lambda_{A,i}\lambda_{B,i}|) - (\lambda_{A,i}^2 - |\lambda_{A,i}\lambda_{B,i}|) = |\lambda_{B,i}| \cdot ( |\lambda_{B,i}| - |\lambda_{A,i}|)<0,\]
and hence
\[\mu_{\min} \ge \min_{i} \{ \lambda_{A,i}^2 + \lambda_{B,i}^2 - 2|\lambda_{A,i}\lambda_{B,i}| \} \ge 64 \tau^2 \]
since $\tau \le 1/8$.

For $d$ dimensions we can carry out the similar argument, except that $\lambda_{A,i}$ and $\lambda_{B,i}$ are replaced by the corresponding results of the new matrix. By the properties of tensor products, the eigenvalues of $\bb{L}_{h,d}$ can be represented by the sum of the eigenvalues of $\bb{L}_h$ as
\[\lambda_{h,d} = \lambda_{h,i_1} + \lambda_{h,i_2} + \cdots + \lambda_{h,i_d}, \quad  \lambda_{h,i_j} \text{~is the eigenvalue of~$\bb{L}_h$}. \]
Hence the eigenvalues with the smallest and largest absolute values are
\[\lambda_{h,1}^d = -4d \sin^2 \frac{\pi}{2N_x}, \qquad \lambda_{h,N_x-1}^d = -4d \sin^2 \frac{(N_x-1)\pi}{2N_x}.\]
In this case, one only requires that the time step satisfy $\tau \le 1/(8d)$ for both cases, and $\beta \le 1/(4d)$ for the explicit scheme, which therefore implies the desired estimates.

(2) The classical method is to iteratively solve
$\bb{A}\bb{u}^{n+1} = \bb{B}\bb{u}^n + \bb{f}^{n+1}$,
where $\bb{A}$ is of order $N_A = (N_x-1)^d$. When $\theta =0, 1/2$, the matrix $\bb{B}$ has sparsity number $s_B\sim d$. The number of basic operations involved on the right-hand side is $\mathcal{O}(N_As_B) = \mathcal{O}(dN_x^d)$. For $\theta = 0$, the total run time is
\[C_{\text{Euler}} = N_t \cdot \mathcal{O}(dN_x^d) = \mathcal{O}(d N_tN_x^d) = \mathcal{O}(d^2 N_x^{d+2}), \quad \theta = 0.\]
For $\theta = 1/2$, noting that $\bb{A}$ is positive definite with sparsity number $s_A \sim d$ and condition number
\[\kappa_A = \frac{\lambda_{A, \max}}{\lambda_{A, \min}}
\sim \frac{1+4d\theta\beta}{1+8d\theta \beta h^2}
\sim \frac{1+4d\theta\beta}{1+8d\theta \tau } = \mathcal{O}(h^{-1}),\]
hence the time complexity of the CG method is
\[\mathcal{O}(N_As_A\sqrt{\kappa_A}) = \mathcal{O}(dN_x^{d+0.5}).\]
Thus the total run time is
\[C_{\text{CN}}  = N_t \mathcal{O}(d N_x^d + dN_x^{d+0.5})) = \mathcal{O}( d N_t N_x^{d+0.5}))
        = \mathcal{O}( d^2 N_x^{d+1.5} ), \quad \theta = 1/2.\]

(3) According to the established result in (1), the condition number of $\bb{H}$ is
\[\kappa \le \begin{cases}
\frac{3}{8 d \tau} \sim  N_x^2, \quad & \theta = 0, \\
\frac{4d \beta}{8 d \tau} \sim  N_x^2, \quad & \theta = 1/2.
\end{cases}\]
One can check that the sparsity number $s \sim d$. The time complexity is
\[Q  = \mathcal{O} ( d N_x^2 \log(N_x^2/d)  ) .\]
This completes the proof.
\end{proof}

It is observed that the time complexity of the classical algorithm depends on $N_x^d$, where $d$ is the spatial dimension, while the quantum algorithm only depends on $N_x$. This shows that the quantum difference method has exponential acceleration with respect to the spatial dimension $d$. On the other hand, unlike the classical difference methods, the quantum treatment of the heat equation gives the same time complexity for the forward Euler and Crank-Nicolson discretizations, which is not reported in the literature.

\begin{remark}
We cannot derive the spectral norm dependence on the coefficient matrix in the condition number from Theorem 7 of \cite{Berry-2014}. For simplicity we consider the two-dimensional case. The resulting ODEs are given in \eqref{odeheat}, where $A = \frac{1}{h^2} \bb{L}_h$ corresponds to the one in \cite{Berry-2014}. For this specific matrix, the condition for the time step in \cite[Theorem 7]{Berry-2014} is $\tau = 1/\|A\| \sim h^2$, which simply refers to the step size of the explicit scheme.
\end{remark}

\section{First order hyperbolic equation} \label{sec:1ohyperblic}

For hyperbolic problems, a digital quantum algorithm combining the finite volume method and the reservoir
technique for symmetric first-order linear hyperbolic systems can be found in \cite{Fillion-Lorin-2019}. In this section, we consider the quantum difference methods for solving the first order hyperbolic equation
\[u_t + u_{x_1} + u_{x_2} + \cdots + u_{x_d} = 0, \quad \bb{x} = (x_1,x_2,\cdots,x_d) \in (0,1)^d\]
in $d$ dimensions with homogeneous inflow boundary conditions, where $u = u(t,x_1,x_2,\cdots, x_d)$.

Let $\bb{j} = (j_1,j_2,\cdots, j_d)$. The upwind scheme can be written as
\[\frac{u_{\bb{j}}^{n+1}-u_{\bb{j}}^n}{\tau} + \sum\limits_{k=1}^d \frac{u^n_{\bb{j}}- u^n_{\bb{j}-\bb{e}_k}}{h} = 0\]
or
\[u_{\bb{j}}^{n+1} + (d\beta-1)u_{\bb{j}}^n - \beta \sum\limits_{k=1}^d u^n_{\bb{j}-\bb{e}_k} = 0,\quad \beta = \tau/h,\]
where $\bb{e}_k = (0,\cdots,0,1,0,\cdots,0)$ with $k$-th entry being 1, $j_k = 1,2,\cdots, N_x$ and $n = 0,1,\cdots,N_t-1$. One easily finds that the above system can be written in matrix form as
\[\bb{u}^{n+1} - \bb{B} \bb{u}^n = \bb{0},\]
with
\begin{align*}
\bb{B}
& =  \lambda(\bb{I} \otimes \bb{I} \otimes \cdots \otimes \bb{T}_h
+ \bb{I} \otimes \cdots \otimes \bb{T}_h \otimes \bb{I} + \cdots + \bb{T}_h\otimes \bb{I} \otimes \cdots \otimes \bb{I}) \\
& \quad  + (1-d\lambda) \underbrace{ \bb{I}\otimes\bb{I} \otimes \cdots  \otimes \bb{I}}_{d~\text{matrices}},
 \end{align*}
where
\[\bb{T}_h = \begin{bmatrix}
0   &         &            &         &    \\
1   &    0    &            &         &    \\
    & \ddots  &  \ddots    &         &    \\
    &         &  \ddots    &  \ddots &    \\
    &         &            &   1     &  0  \\
\end{bmatrix}_{N_x \times N_x},\]
and $\bb{u}$ represents the vector form of the $d$-order tensor $\bb{U} = (u_{\bb{j}}) = (u_{j_1,j_2,\cdots, j_d})$.
The coefficient matrix for the quantum difference method is then given by
\begin{equation}\label{ceffL}
\bb{L} =
\begin{bmatrix}
\bb{I}  &            &           &            \\
-\bb{B} & \bb{I}     &           &            \\
        &\ddots      & \ddots    &    \\
        &            & -\bb{B}   & \bb{I}     \\
\end{bmatrix}.
\end{equation}

To give the bounds of the eigenvalues of the corresponding matrix $\bb{H}$, or equivalently the singular values of $\bb{L}$, we introduce the Gershgorin-type theorem for singular values.
\begin{lemma}\cite{Qi-1984} \label{lem:Qi}
Let $\bb{A} = (a_{ij})$  be a square matrix of order $n$. Write
\[r_i = \sum\limits_{j\ne i} |a_{ij}|, \quad c_i = \sum\limits_{j\ne i} |a_{ji}|, \quad  s_i = \max (r_i, c_i)\]
for $i=1,2,\cdots, n$. Then each singular value of $\bb{A}$ lies in one of the real intervals
\[[(|a_{ii}|-s_i)_+,  |a_{ii}|+s_i],\]
where $a_+ = \max (0,a)$.
\end{lemma}

The error of the upwind scheme is $\mathcal{O}(\tau + d h)$ if the exact solution is regular. In the following theorem we set $\delta = \mathcal{O}(dh) = \mathcal{O}(d/N_x)$.

\begin{theorem}\label{thm:hyperbolic}
Let $\beta = \tau/h \le 1/d$. Then the condition number and the sparsity of $\bb{H}$ satisfy
\[\kappa = \mathcal{O}(N_t) \qquad \mbox{and} \qquad s = \mathcal{O}(d).\]
For fixed spatial step $h$, let $\tau = h/d$. Then the time complexities of the classical difference methods and the quantum difference method for solving the $d$-dimensional first order hyperbolic equation are
\[C = \mathcal{O}(d^2 N_x^{d+1})  \quad \text{and} \quad
Q = \mathcal{O}(d^2N_x \log(N_x/d) ),\]
respectively.
\end{theorem}
\begin{proof}
(1) For ease of presentation, we only consider the case of $d=3$. In this case, $\bb{B}$ has the following form
\[ \bb{B}  = {\scriptsize \left[\begin{array}{cccc:cccccccc}
\tilde{\bb{T}}_h    &          &              &         &       &        &              &  \\
\lambda\bb{T}_h  & \tilde{\bb{T}}_h   &              &         &            &   &              &   \\
          & \ddots   & \ddots       &         &            &        &         & \\
          &          & \lambda\bb{T}_h    & \tilde{\bb{T}}_h  &            &        &              &  \\
 \hdashline
\lambda\bb{I}     &         &          &        & \tilde{\bb{T}}_h    &          &          & \\
           & \lambda\bb{I}  &          &        &  \lambda\bb{T}_h   & \tilde{\bb{T}}_h   &          & \\
           &         & \ddots   &        &             &   \ddots       & \ddots   & \\
           &         &          & \lambda\bb{I} &             &          &    \lambda\bb{T}_h      &  \tilde{\bb{T}}_h \\
\end{array} \right]}, \qquad \tilde{\bb{T}}_h = (1-3\lambda)\bb{I}+\lambda\bb{T}_h,\]
where the repeated blocks are omitted. Applying Lemma \ref{lem:Qi} to get
\[\|\bb{B}\| = \sigma_{\max}(\bb{B}) \le 1, \qquad \sigma_{\max}(\bb{L}) \le 2.\]

By definition, $\sigma_{\min}(\bb{L}) = 1/\sigma_{\max}(\bb{L}^{-1})$. After simple algebra, one has
\[
{ \scriptsize
\bb{L}^{-1} =
\begin{bmatrix}
\bb{I}          &              &                 &          &            \\
 \bb{B}        & \bb{I}       &                 &          &        \\
 \bb{B}^2      &  \ddots      &  \ddots         &          &\\
   \vdots       &  \ddots      &  \ddots         &  \ddots  &   \\
\bb{B}^{N_t-1} &   \cdots     &   \bb{B}^2     &  \bb{B} &  \bb{I}      \\
\end{bmatrix}
= \begin{bmatrix}
\bb{I}          &              &                 &          &            \\
                & \bb{I}       &                 &          &        \\
                &              &  \ddots         &          &\\
                &              &                 &  \ddots  &   \\
                &              &                 &          &  \bb{I}      \\
\end{bmatrix}+
\begin{bmatrix}
                 &              &                 &          &            \\
 \bb{B}        &              &                 &          &        \\
                 &  \ddots      &                 &          &\\
                 &              &  \ddots         &          &   \\
                 &              &                 &  \bb{B} &         \\
\end{bmatrix} + \cdots,}
\]
which gives
\begin{align*}
\sigma_{\max}(\bb{L}^{-1})
& = \|\bb{L}^{-1}\| \le \|\bb{I}\| + \|\bb{B}\| + \cdots +\|\bb{B}^{N_t-1}\| \\
& \le \|\bb{I}\| + \|\bb{B}\| + \|\bb{B}\|^2 +\cdots +\|\bb{B}\|^{N_t-1} \le N_t = 1/\tau.
\end{align*}
Thus, $\sigma_{\min}(\bb{L}) \ge \tau$, as required.

(2) The classical method is to iteratively solve $\bb{u}^{n+1} = \bb{B}\bb{u}^n$,
where $\bb{B}$ is of order $n_B = N_x^d$. One can check that  matrix $\bb{B}$ has sparsity number $s_B\sim d$.
Thus, the number of basic operations involved on the right-hand side is $\mathcal{O}(n_Bs_B) = \mathcal{O}(dN_x^d)$, and the total run time is
\[C = N_t \cdot \mathcal{O}(dN_x^d) = \mathcal{O}(d N_tN_x^d) = \mathcal{O}(d^2 N_x^{d+1}).\]

(3) The condition number is
\[\kappa \le 2/\tau \sim N_t = dN_x.\]
One can check that the sparsity number $s \sim d$, then the query complexity is
\[Q = \mathcal{O}(d^2N_x \log(N_x/d) ).\]
This completes the proof.
\end{proof}

\section{Multiscale problem} \label{sec:multiscale}

Consider the multiscale hyperbolic heat system (telegraph equations)
\begin{equation} \label{heatstiff}
\begin{cases}
u_t + v_x = 0, \\
v_t + \frac{1}{\varepsilon} u_x = - \frac{1}{\varepsilon}v, \quad 0<\varepsilon \ll 1,
\end{cases}
\end{equation}
where $\varepsilon$ is the relaxation time or the scaling parameter. For simplicity, we assume the exact solution is regular. Note that numerically approximating the above system is challenging due to the stiffness of the problem for both the convection and collision terms \cite{Jin-Pareschi-Toscani-1998,Jin-Pareschi-Toscani-2000}. The key idea to tackle  problem \eqref{heatstiff} is to reformulate it as a (nonstiff) linear hyperbolic system with stiff relaxation term called the diffusive relaxation system \cite{Jin-Pareschi-Toscani-1998}:
\begin{align}\label{relaxation}
\begin{cases}
u_t + v_x = 0, \\
v_t + u_x = - \frac{1}{\varepsilon}\Big( v + (1-\varepsilon) u_x \Big), \quad 0<\varepsilon \ll 1.
\end{cases}
\end{align}
This system \eqref{relaxation} has the form of Jin-Xin relaxation model used to construct Riemann solver free shock capturing schemes for conservation laws by Jin and Xin (cf. \cite{Jin-Xin-1995}).

\subsection{AP schemes}

The AP schemes have been developed for a wide range of time-dependent kinetic and hyperbolic equations. The fundamental idea is to design numerical methods that preserve the asymptotic limits from the microscopic to the macroscopic models in the discrete setting \cite{JinReview-2012, JinActa}. We first consider two typical AP schemes presented in \cite{JinReview-2012,Hu-Jin-Li-2017} for solving \eqref{heatstiff} or \eqref{relaxation}.

\subsubsection{IMEX scheme}

The first AP scheme is the following implicit-explicit (IMEX) scheme (cf. (2.15) in \cite{JinReview-2012}):
\begin{equation} \label{imex}
\begin{cases}
\frac{u_j^{n+1} - u_j^n}{\tau} + \frac{v_{j+1}^n-v_{j-1}^n}{2h} - \frac{h}{2}\frac{u_{j-1}^n-2u_j^n+u_{j+1}^n}{h^2}= 0, \\
\frac{v_j^{n+1} - v_j^n}{\tau} + \frac{u_{j+1}^n-u_{j-1}^n}{2h} - \frac{h}{2}\frac{v_{j-1}^n-2v_j^n+v_{j+1}^n}{h^2}
    = - \frac{1}{\varepsilon}\Big( v_j^{n+1} + (1-\varepsilon)\frac{u_{j+1}^{n+1}-u_{j-1}^{n+1}}{2h} \Big),
\end{cases}
\end{equation}
where $j = 1,\cdots,N_x-1$ and $n = 0,\cdots, N_t-1$. Notice that the scheme treats the relaxation term implicitly. On a classical computer, one can obtain $u^{n+1}$ from the first equation of \eqref{imex} explicitly and then substitute it into the right-hand side of the second equation. The remaining equation for $v^{n+1}$ can also be implemented explicitly.

In order to write it as a large linear system as in \eqref{ODEforwardAxb}, define
\[\bb{u}(t) = [u_1(t), \cdots, u_{N_x-1}(t)]^T, \quad  \bb{v}(t) = [v_1(t), \cdots, v_{N_x-1}(t)]^T.\]
The corresponding spatial discretization of \eqref{imex} is
\begin{equation*}
\begin{cases}
\frac{{\rm d}}{{\rm d}t} \bb{u}(t)  + \frac{1}{2h} \bb{M}_h \bb{v}(t) - \frac{1}{2h} \bb{L}_h \bb{u}(t) - \frac{1}{2h} (\bb{b}(t) - \widetilde{\bb{c}}(t)) = \bb{0}, \\
\frac{{\rm d}}{{\rm d}t} \bb{v}(t)  + \frac{1}{2h} \bb{M}_h \bb{u}(t) - \frac{1}{2h} \bb{L}_h \bb{v}(t) + \frac{1}{2h} (\widetilde{\bb{b}}(t) - \bb{c}(t))
   =-\frac{1}{\varepsilon}\Big( \bb{v}(t) + \frac{1-\varepsilon}{2h}\bb{M}_h\bb{u}(t) + \frac{1-\varepsilon}{2h}\widetilde{\bb{b}}(t) \Big),
\end{cases}
\end{equation*}
where,
\begin{equation}\label{MhLh}
\bb{M}_h =
\begin{bmatrix}
0  &  1       &           &      &    \\
 -1  & 0       & \ddots    &      &    \\
    &  \ddots  & \ddots    &  \ddots    &    \\
    &          & \ddots    & 0   & 1  \\
    &          &           &  -1   & 0 \\
\end{bmatrix}
, \qquad
\bb{L}_h =
\begin{bmatrix}
-2  &  1       &           &      &    \\
 1  & -2       & \ddots    &      &    \\
    &  \ddots  & \ddots    &  \ddots    &    \\
    &          & \ddots    & -2   & 1  \\
    &          &           &  1   & -2 \\
\end{bmatrix}_{(N_x-1) \times (N_x-1)},
\end{equation}
\[
\bb{b}(t) =
\begin{bmatrix}
u_0(t) \\
0  \\
\vdots\\
0 \\
u_{N_x}(t) \\
\end{bmatrix}, \qquad
\widetilde{\bb{b}}(t) =
\begin{bmatrix}
-u_0(t) \\
0  \\
\vdots\\
0 \\
u_{N_x}(t) \\
\end{bmatrix}, \qquad
\bb{c}(t) =
\begin{bmatrix}
v_0(t) \\
0  \\
\vdots\\
0 \\
v_{N_x}(t) \\
\end{bmatrix}, \qquad
\widetilde{\bb{c}}(t) =
\begin{bmatrix}
-v_0(t) \\
0  \\
\vdots\\
0 \\
v_{N_x}(t) \\
\end{bmatrix}.
\]
Further time approximation yields
\begin{equation*}
\begin{cases}
\frac{\bb{u}^{n+1} - \bb{u}^n}{\tau}  + \frac{1}{2h} \bb{M}_h \bb{v}^n - \frac{1}{2h} \bb{L}_h \bb{u}^n - \frac{1}{2h} (\bb{b}^n - \widetilde{\bb{c}}^n) = \bb{0}, \\
\frac{\bb{v}^{n+1} - \bb{v}^n}{\tau}  + \frac{1}{2h} \bb{M}_h \bb{u}^n - \frac{1}{2h} \bb{L}_h \bb{v}^n + \frac{1}{2h} (\widetilde{\bb{b}}^n - \bb{c}^n)
   =-\frac{1}{\varepsilon}\Big( \bb{v}^{n+1} + \frac{1-\varepsilon}{2h}\bb{M}_h\bb{u}^{n+1} + \frac{1-\varepsilon}{2h}\widetilde{\bb{b}}^{n+1} \Big).
\end{cases}
\end{equation*}
Let $\beta = \tau/h$. The above system can be written as
\begin{equation}\label{imex0}
\begin{cases}
-\bb{B}\bb{u}^n + \bb{u}^{n+1} + \bb{A}\bb{v}^n = \bb{f}^{n+1}, \\
-\bb{B}\bb{v}^n + \gamma \bb{v}^{n+1} + \bb{A}\bb{u}^n + \nu \bb{A}\bb{u}^{n+1} = \bb{g}^{n+1},
\end{cases}
\end{equation}
where, $\gamma = 1+\tau/\varepsilon$, $\nu = (1-\varepsilon)/\varepsilon$, and
\[\bb{A} = \frac{\beta}{2}\bb{M}_h, \qquad \bb{B} = \bb{I} + \frac{\beta}{2} \bb{L_h},\]
\[\bb{f}^{n+1} = \frac{\beta}{2}(\bb{b}^n-\widetilde{\bb{c}}^n), \qquad \bb{g}^{n+1} = \frac{\beta}{2}(\bb{c}^n - \widetilde{\bb{b}}^n - \nu \widetilde{\bb{b}}^{n+1}). \]

Introducing the following notations
\[\bb{U} = [\bb{u}^1; \cdots ;\bb{u}^{N_t}], \qquad \bb{V} = [\bb{v}^1; \cdots ;\bb{v}^{N_t}], \qquad \bb{S} = [\bb{U}; \bb{V}],\]
one then gets a linear system
\begin{equation}\label{APimex}
\bb{L}_{\text{IMEX}} \bb{S} = \bb{F}_{\text{IMEX}},
\end{equation}
where, $\bb{L}_{\text{IMEX}} = (\bb{L}_{ij})_{2\times 2}$, $\bb{F}_{\text{IMEX}} = [\bb{F}_1; \bb{F}_2]$, and
\[\bb{L}_{11} =
\begin{bmatrix}
\bb{I}  &            &           &            \\
-\bb{B} & \bb{I}     &           &            \\
        &\ddots      & \ddots    &    \\
        &            & -\bb{B}   & \bb{I}     \\
\end{bmatrix}
, \qquad
\bb{L}_{12} =
\begin{bmatrix}
\bb{O}  &            &           &            \\
\bb{A} & \bb{O}     &           &            \\
        &\ddots      & \ddots    &    \\
        &            & \bb{A}   & \bb{O}     \\
\end{bmatrix}, \qquad
\bb{F}_1 =
\begin{bmatrix}
\bb{f}^1 + \bb{B}\bb{u}^0-\bb{A}\bb{v}^0 \\
\bb{f}^2  \\
\vdots\\
\bb{f}^{N_t} \\
\end{bmatrix}
\]
\[\bb{L}_{21} =
\begin{bmatrix}
\nu\bb{A}  &            &           &            \\
\bb{A} & \nu\bb{A}     &           &            \\
        &\ddots      & \ddots    &    \\
        &            & \bb{A}   & \nu\bb{A}     \\
\end{bmatrix}
, \qquad
\bb{L}_{22} =
\begin{bmatrix}
\gamma\bb{I}  &            &           &            \\
-\bb{B} & \gamma\bb{I}     &           &            \\
        &\ddots      & \ddots    &    \\
        &            & -\bb{B}   & \gamma\bb{I}    \\
\end{bmatrix}
, \qquad
\bb{F}_2 =
\begin{bmatrix}
\bb{g}^1 - \bb{A}\bb{u}^0 + \bb{B}\bb{v}^0\\
\bb{g}^2  \\
\vdots\\
\bb{g}^{N_t} \\
\end{bmatrix}.
\]

On a classical computer, one can obtain good results by solving the IMEX scheme directly for a relatively large -- compared with $\varepsilon$ -- time steps and mesh sizes. However, the classical simulation becomes difficult when one solves \eqref{APimex} with the QLSA, because the condition number of the matrix $\bb{H}$ is very large. In fact, for fixed step sizes $\tau$ and $h$, the eigenvalues of $\bb{H}$ corresponding to \eqref{APimex} satisfy
\[|\lambda|_{\max} \sim \frac{1}{\varepsilon}, \qquad |\lambda|_{\min} \le 1.\]
A simple argument is described as follows.  By the characteristics of the matrix, it is sufficient to consider only the computation from  $t_n$ to $t_{n+1}$. The coefficient matrix in this case is
\[\bb{L} = \begin{bmatrix}  \bb{I} & \bb{O} \\ \nu \bb{A} & \gamma \bb{I}, \end{bmatrix}, \]
and
\[\bb{D} := \bb{L}\bb{L}^T = \begin{bmatrix}   \bb{I}  & -\nu\beta/2\bb{M}_h \\  \nu\beta/2\bb{M}_h   & \gamma^2\bb{I}-\nu^2\beta^2/4\bb{M}_h^2 \end{bmatrix}.\]
The matrix $\bb{M }_h$ is a standard antisymmetric matrix whose eigenvalue $\lambda_h$ is zero or pure imaginary, and thus $\lambda_h^2 \le 0$. By $\|\bb{M}_h\|_\infty = 2$, $|\lambda_h| \le 2$. Let $\Lambda$ be the matrix consisting of these eigenvalues. One has
\[\bb{D} \sim  \begin{bmatrix}   \bb{I}  & -\nu\beta/2\Lambda \\  \nu\beta/2\Lambda   & \gamma^2\bb{I}-\nu^2\beta^2/4\Lambda^2 \end{bmatrix},\]
which implies that the maximum eigenvalue of $\bb{D}$ satisfies
\[\mu_{\max} \le \|\bb{D}\|_\infty \le \frac{\nu^2\beta^2}{4} \cdot |\lambda_h|_{\max}^2 + \gamma^2 + \frac{\nu\beta}{2} |\lambda_h|_{\max} \le C(\tau,h) \frac{1}{\varepsilon^2}. \]
By the Rayleigh quotient theorem for symmetric matrices \cite{Horn2013},
\[\mu_{\max} \ge \max_{i} \bb{D}_{ii} \ge C(\tau,h) \frac{1}{\varepsilon^2}, \qquad \mu_{\min} \le \min_{i} \bb{D}_{ii} \le 1.\]
The claim follows from the above two equations.

A simple way to tackle  this problem is to apply a preconditioner, that is, we can find a simple invertible matrix $\bb{P}$ such that $(\bb{P}\bb{L}_{\text{IMEX}})\bb{S} = \bb{P} \bb{F}_{\text{IMEX}}$ can be solved efficiently. \\

We remark that in the presence of the preconditioner $\bb{P}$, we now assume access to the unitary preparing the initial quantum state proportional to $\bb{P}\bb{F}_{\text{IMEX}}$, instead of $\bb{F}_{\text{IMEX}}$. This assumption is motivated straightforwardly due to the simplicity of $\bb{P}$. For instance, following state preparation protocols like \cite{gleinig2021efficient} or \cite{zhang2022quantum}, the cost only depends on the sparsity and dimension of the states, which are identical for $\bb{F}_{\text{IMEX}}$ and $\bb{P}\bb{F}_{\text{IMEX}}$. Similarly, we assume access to a sparse-access query model for $(\bb{P}\bb{L}_{\text{IMEX}})$ instead of $\bb{L}_{\text{IMEX}}$, where we easily see from $\bb{P}$ that the sparsities of these matrices are identical.

\begin{theorem}\label{thm:precondition}
Let $\bb{P}$ be defined in block form by
\[\bb{P} = \begin{bmatrix} \bb{I} & \bb{O} \\ \bb{O} & \varrho\bb{I} \end{bmatrix},
\qquad \varrho = \frac{\varepsilon}{1+\varepsilon},\]
where $\varrho$ is the preconditioned parameter. Let $\widetilde{\bb{L}}_{\text{IMEX}} = \bb{P}\bb{L}_{\text{IMEX}}$ and $\widetilde{\bb{F}}_{\text{IMEX}} = \bb{P} \bb{F}_{\text{IMEX}}$. Consider solving the following preconditioned linear system
\begin{equation}\label{AP1Pre}
\widetilde{\bb{L}}_{\text{IMEX}}\bb{S} = \widetilde{\bb{F}}_{\text{IMEX}}.
\end{equation}
Let $\lambda$ be the eigenvalue of the Hermitian $\bb{H}$ corresponding to $\widetilde{\bb{L}}_{\text{IMEX}}$, i.e., $\bb{H}=\begin{bmatrix} 0 & \widetilde{\bb{L}}_{\text{IMEX}} \\ \widetilde{\bb{L}}_{\text{IMEX}}^T & 0 \end{bmatrix}$. Then for fixed  $\tau$ and $h$, the upper and lower bounds of $|\lambda|$ are independent of $\varepsilon$.
\end{theorem}
\begin{proof}For convenience, we omit the subscript $\text{IMEX}$ in the following.
The parameters associated with $\varepsilon$ in the matrix $\widetilde{\bb{L}}$ are
\[\varrho = \frac{\varepsilon}{1+\varepsilon},
\quad \varrho \nu = \frac{1-\varepsilon}{1+\varepsilon},
\quad \varrho \gamma = \frac{\tau+\varepsilon}{1+\varepsilon}.
\]
They are all bounded as $\varepsilon \to 0$  (they only depend on the step sizes). One easily finds that $\widetilde{\bb{L}}$ converges to a fixed invertible matrix as $\varepsilon \to 0$. Since the singular values of $\widetilde{\bb{L}}$ are continuous with respect to $\varepsilon$, the upper and lower bounds of the singular values of  matrix $\widetilde{\bb{L}}$ are independent of $\varepsilon$ for fixed $\tau$ and $h$. Moreover, the absolute values of $\lambda$ are exactly the singular values of $\widetilde{\bb{L}}$, thus the upper and lower bounds of $|\lambda|$ are also independent of $\varepsilon$.
\end{proof}

In the simulation, the initial-boundary values of \eqref{heatstiff} are chosen in such a way that the exact solution is
\[u = {\rm e}^{at} \sin (ax),  \quad v = {\rm e}^{at} \cos (ax), \quad a = -1/(1+\varepsilon).\]
The spatial and temporal domains are taken as $[-1,1]$ and $[0,1]$, respectively. We still consider the original HHL algorithm in view of various detailed implementations in the literature, see \cite{CaoCircuit-2012,CaoPoisson-2013} for example.
Given a Hermitian matrix $\bb{H}$, let $\bb{U} = \mathrm{e}^{\mathrm{i} \bb{H} t_0}$ be the unitary matrix in the quantum phase estimation. Denote the eigenvalue of $\bb{H}$ by $\lambda$. Then the corresponding phase $\phi \in [0,1)$ is defined by
${\rm e}^{2\pi{\rm i} \phi} = {\rm e}^{{\rm i} \lambda t_0}$, where ${\rm i} = \sqrt{-1}$ and $t_0$ is the evolution time. Noting that $\ln(\mathrm{e}^{{{\rm i} \theta}})/{\rm i} = \theta$ for $\theta \in [-\pi, \pi]$, we require that $|\lambda| t_0 \le \pi$. Then the eigenvalue $\lambda$ can be represented by the phase $\phi$ as \[\lambda = \begin{cases}
\frac{2\pi \phi}{t_0}, \quad & \lambda > 0, \\
\frac{2\pi(\phi - 1)}{t_0}, \quad & \lambda < 0.
\end{cases}\]
The evolution time $t_0$ will be set as $t_0 = \frac{2\pi}{2^{n_t}}C_t$, where $n_t$ is the number of qubits in the clock register, and $C_t = 10^p$ with $p = [\log_{10} ( {2^{n_t-1}}/{|\lambda|_{\max}}) ]$ being an integer. The reason for this choice is that the integer part of $\lambda C_t$ can be represented exactly by $n_t$-bits, while $C_t$ shifts the fractional part of $\lambda$ to the integer part.

\begin{figure}[!ht]
  \centering
  \subfigure[$\varepsilon = 10^{-1}$]{\includegraphics[scale=0.35]{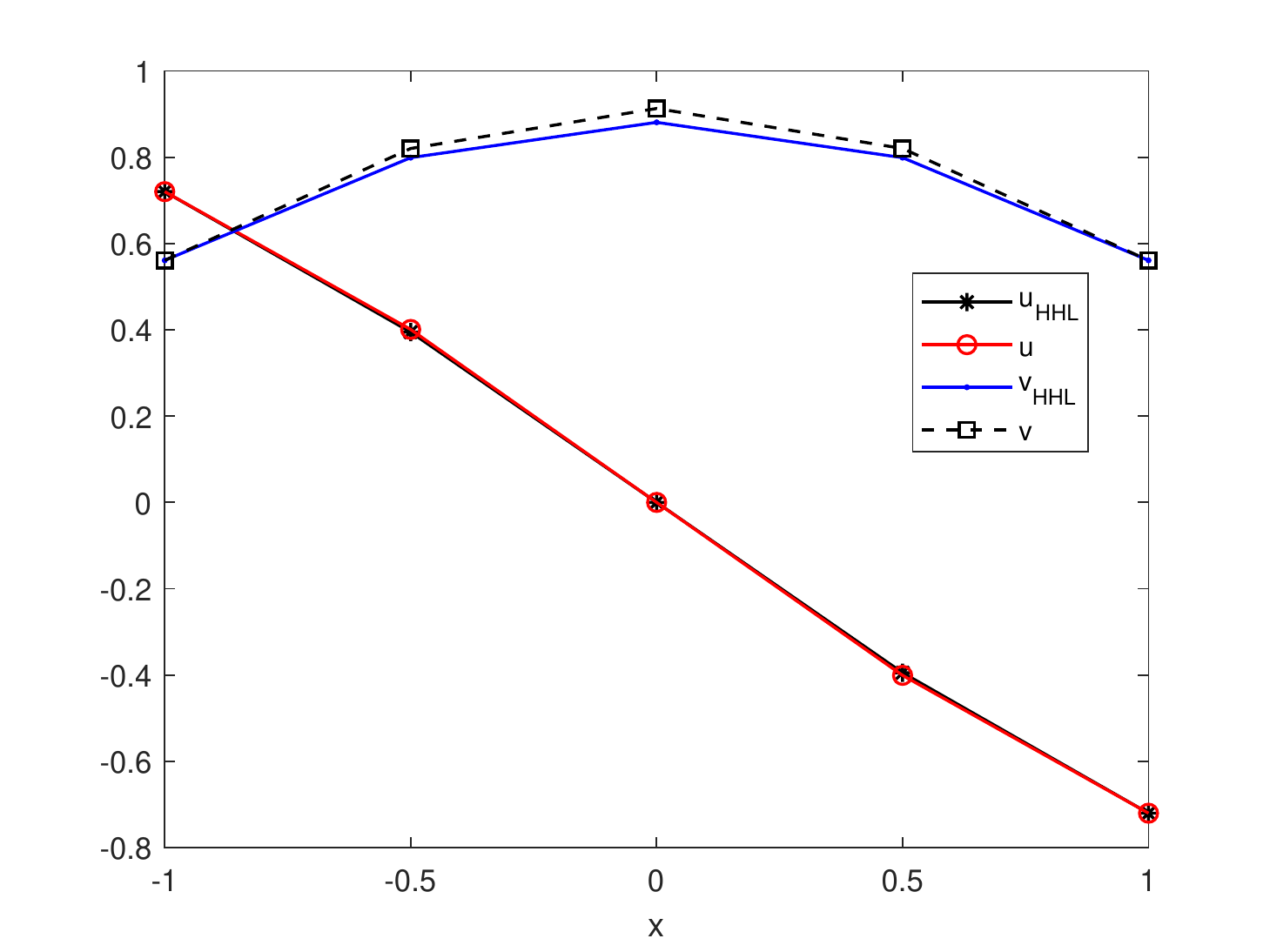}}
  \subfigure[$\varepsilon = 10^{-4}$]{\includegraphics[scale=0.35]{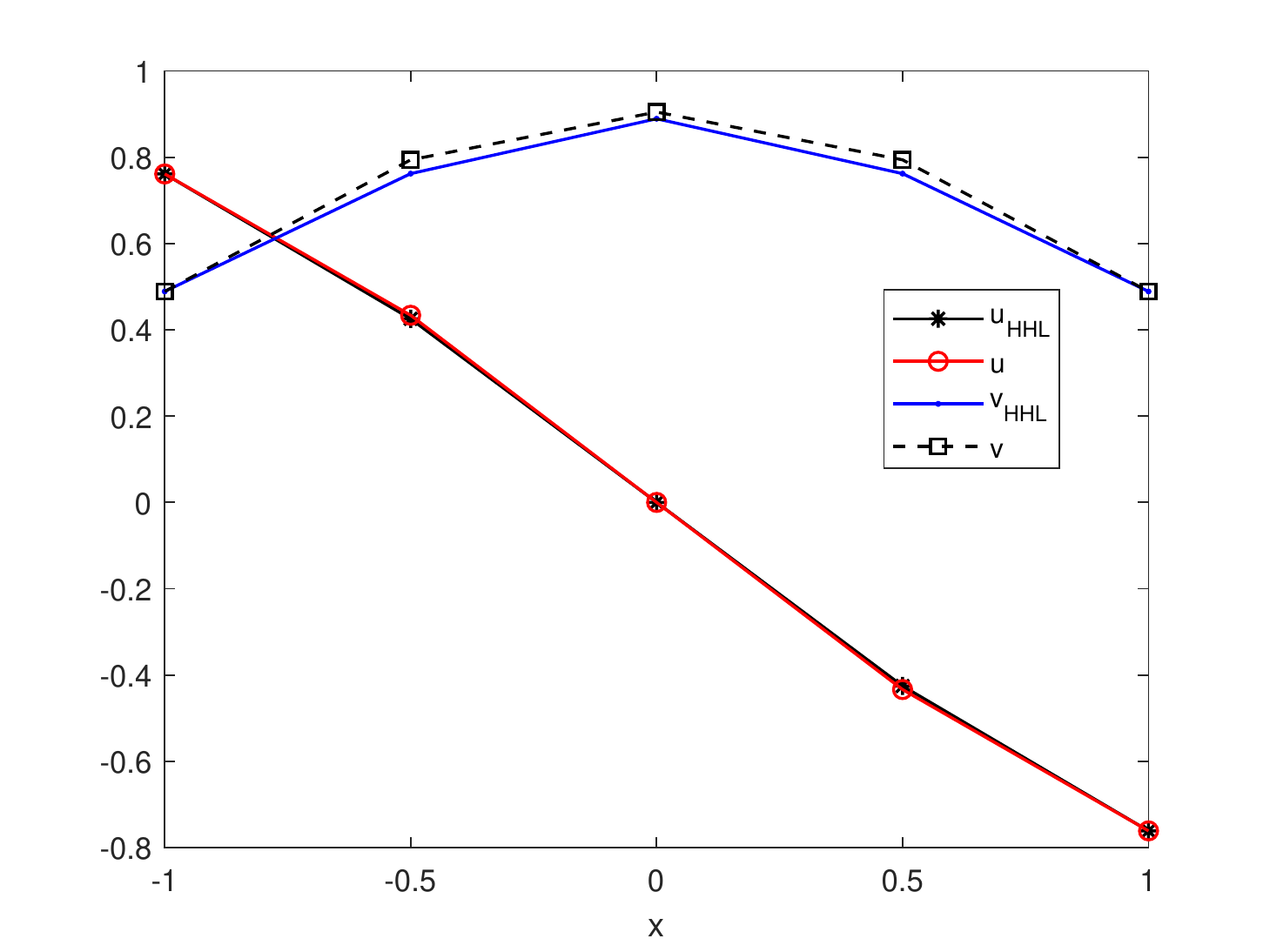}}
  \subfigure[$\varepsilon = 10^{-8}$]{\includegraphics[scale=0.35]{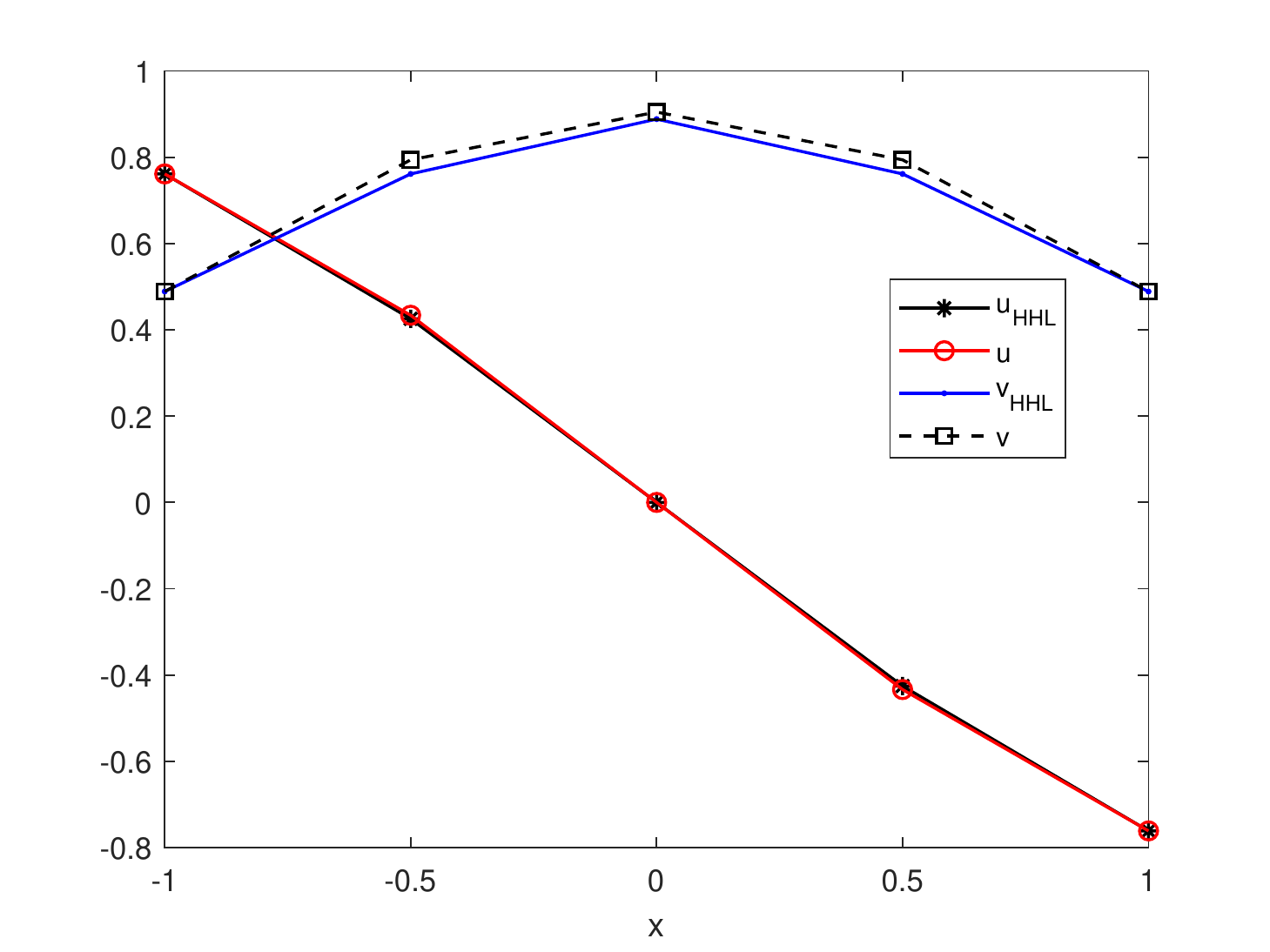}}\\
  \caption{Exact and HHL solutions of $u$ and $v$ for the preconditioned IMEX scheme with $\tau=0.05$, $h = 0.5$ and $\varrho = \varepsilon/(1+\varepsilon)$}\label{fig:ApheatPre}
\end{figure}

The number of qubits in the clock register is $n_t = 10$. For $\varepsilon = 10^{-i}, i = 1,4, 8$, Fig.~\ref{fig:ApheatPre} displays the numerical and exact solutions of $u$ and $v$, from which one observes that the preconditioned scheme is $\varepsilon$-independent, which verifies the conclusion in Theorem \ref{thm:precondition}.

Unless otherwise specified, the notations in this subsection will be frequently used in the following text.

\subsubsection{Diffusive relaxation scheme}

The second AP scheme is the diffusive relaxation scheme proposed by Jin, Pareschi and Toscani in \cite{Jin-Pareschi-Toscani-1998} for solving multiscale discrete-velocity kinetic equations.

The diffusive relaxation scheme has two steps:
\begin{enumerate}
  \item Relaxation step
  \begin{equation} \label{relaxationstep}
\begin{cases}
\frac{u_j^{*} - u_j^n}{\tau} = 0, \quad j = 0,\cdots, N_x, \\
\frac{v_j^{*} - v_j^n}{\tau}  = - \frac{1}{\varepsilon}\Big( v_j^{*} + (1-\varepsilon)\frac{u_{j+1}^{*}-u_{j-1}^{*}}{2h} \Big),
\quad j = 1,\cdots, N_x-1.
\end{cases}
\end{equation}
 In the computation, we take $v_0^* = v_0^n$ and $v_{N_x}^* = v_{N_x}^n$.
  \item Convection step
  \begin{equation} \label{convectionstep}
\begin{cases}
\frac{u_j^{n+1} - u_j^*}{\tau} + \frac{v_{j+1}^*-v_{j-1}^*}{2h} - \frac{h}{2}\frac{u_{j-1}^*-2u_j^*+u_{j+1}^*}{h^2}= 0, \\
\frac{v_j^{n+1} - v_j^*}{\tau} + \frac{u_{j+1}^*-u_{j-1}^*}{2h} - \frac{h}{2}\frac{v_{j-1}^*-2v_j^*+v_{j+1}^*}{h^2} = 0.
\end{cases}
\end{equation}
\end{enumerate}

The equations \eqref{relaxationstep} and \eqref{convectionstep} can be written in vector form as
\[
\begin{cases}
 \bb{u}^* = \bb{u}^n, \\
 \gamma\bb{v}^* = \bb{v}^n - \nu \bb{A}\bb{u}^n - \nu \frac{\beta}{2} \widetilde{\bb{b}}^n
\end{cases}
\]
and
\[
\begin{cases}
\bb{u}^{n+1} = \bb{B} \bb{u}^* - \bb{A}\bb{v}^* + \frac{\beta}{2}(\bb{b}^n - \widetilde{\bb{c}}^n) , \\
\bb{v}^{n+1} = \bb{B} \bb{v}^* - \bb{A}\bb{u}^* + \frac{\beta}{2}(\bb{c}^n - \widetilde{\bb{b}}^n),
\end{cases}
\]
respectively. Introducing the following notations
\begin{equation}\label{notationA1}
\bb{A}_1 = \frac{1}{\gamma}\bb{A}, \quad \bb{A}_2 = \frac{1}{\gamma} \bb{B}, \quad
\bb{B}_1 = \bb{B} + \frac{\nu}{\gamma}\bb{A}^2, \quad
\bb{B}_2 = \bb{A} + \frac{\nu}{\gamma}\bb{BA},
\end{equation}
\[\widetilde{\bb{f}}^{n+1} = \frac{\beta}{2}(\bb{b}^n - \widetilde{\bb{c}}^n) + \frac{\nu\beta}{2\gamma}\bb{A}\widetilde{\bb{b}}^n, \quad
\widetilde{\bb{g}}^{n+1} = \frac{\beta}{2}(\bb{c}^n - \widetilde{\bb{b}}^n) - \frac{\nu\beta}{2\gamma}\bb{B}\widetilde{\bb{b}}^n,\]
and eliminating $\bb{u}^*$ and $\bb{v}^*$, one obtains the linear system
\begin{equation}\label{APsplit}
\bb{L}_{\text{relaxation}} \bb{S} = \bb{F}_{\text{relaxation}},
\end{equation}
where, $\bb{L}_{\text{relaxation}} = (\bb{L}_{ij})_{2\times 2}$, $\bb{F}_{\text{relaxation}} = [\bb{F}_1; \bb{F}_2]$, and
\[\bb{L}_{11} =
\begin{bmatrix}
\bb{I}  &            &           &            \\
-\bb{B}_1 & \bb{I}     &           &            \\
        &\ddots      & \ddots    &    \\
        &            & -\bb{B}_1   & \bb{I}     \\
\end{bmatrix}
, \qquad
\bb{L}_{12} =
\begin{bmatrix}
\bb{O}  &            &           &            \\
\bb{A}_1 & \bb{O}     &           &            \\
        &\ddots      & \ddots    &    \\
        &            & \bb{A}_1   & \bb{O}     \\
\end{bmatrix}, \qquad
\bb{F}_1 =
\begin{bmatrix}
\widetilde{\bb{f}}^1 + \bb{B}_1\bb{u}^0-\bb{A}_1\bb{v}^0 \\
\widetilde{\bb{f}}^2  \\
\vdots\\
\widetilde{\bb{f}}^{N_t} \\
\end{bmatrix}
\]
\[\bb{L}_{21} =
\begin{bmatrix}
\bb{O}   &            &            &            \\
\bb{B}_2 & \bb{O}     &            &            \\
         & \ddots     & \ddots     &    \\
         &            & \bb{B}_2   & \bb{O}     \\
\end{bmatrix}
, \qquad
\bb{L}_{22} =
\begin{bmatrix}
\bb{I}  &            &           &            \\
-\bb{A}_2 & \bb{I}     &           &            \\
        &\ddots      & \ddots    &    \\
        &            & -\bb{A}_2   & \bb{I}    \\
\end{bmatrix}
, \qquad
\bb{F}_2 =
\begin{bmatrix}
\widetilde{\bb{g}}^1 - \bb{B}_2\bb{u}^0 + \bb{A}_2\bb{v}^0\\
\widetilde{\bb{g}}^2  \\
\vdots\\
\widetilde{\bb{g}}^{N_t} \\
\end{bmatrix}.
\]

\begin{theorem}
Let $\lambda$ be the eigenvalue of the matrix $\bb{H}$ corresponding to the linear system \eqref{APsplit}. Then for fixed step sizes $\tau$ and $h$, the upper and lower bounds of $|\lambda|$ are independent of $\varepsilon$.
\end{theorem}
\begin{proof}
The proof is similar to that of Theorem \ref{thm:precondition}. So we omit it.
\end{proof}

For the preconditioned IMEX scheme \eqref{AP1Pre} with $\tau \sim h^2$, a direct calculation gives $\sigma_{\max}(\widetilde{\bb{L}}) \sim 1$, where $\widetilde{\bb{L}}$ is the coefficient matrix. Under the same condition, one can find that $\sigma_{\max}(\bb{L}_{\text{relaxation}}) \sim h^{-1}$ for the diffusive relaxation scheme \eqref{APsplit} though the condition number is independent of $\varepsilon$. For this reason, we reformulate it as
\begin{equation}\label{AP2Pre}
\begin{bmatrix}
\bb{L}_{11}  & \tau^{-1} \bb{L}_{12} \\
\tau \bb{L}_{12}  & \bb{L}_{22}
\end{bmatrix}
\begin{bmatrix}
\tau^{-1}\bb{U}  \\
\bb{V}
\end{bmatrix}
= \begin{bmatrix}
\tau^{-1}\bb{F}_1  \\
\bb{F}_2
\end{bmatrix},
\end{equation}
where $\tau$ is the time step.

We note that if we wanted to recover the state $[\bb{U}, \bb{V}]^T$ instead without the $\tau^{-1}$ rescaling factor in $[\tau^{-1}\bb{U}, \bb{V}]^T$, we can instead rewrite Eq.~\eqref{AP2Pre} as
\begin{equation}\label{AP2Pre2}
\bb{L}'\begin{bmatrix}
\bb{U}  \\
\bb{V}
\end{bmatrix} \equiv \begin{bmatrix}
\tau^{-1}\bb{L}_{11}  & \tau^{-1} \bb{L}_{12} \\
 \bb{L}_{12}  & \bb{L}_{22}
\end{bmatrix}
\begin{bmatrix}
\bb{U}  \\
\bb{V}
\end{bmatrix}
= \begin{bmatrix}
\tau^{-1}\bb{F}_1  \\
\bb{F}_2
\end{bmatrix},
\end{equation}
where the sparsity of the new matrix $\bb{L}'$ and the matrix denoted $\tilde{\bb{L}}$ in Eq.~\eqref{AP2Pre} are clearly identical, and their condition numbers are comparable, namely $\kappa(\tilde{\bb{L}}), \kappa(\bb{L}') \lesssim \tau^{-1}$ for small $\varepsilon$. This means that the query and gate complexities of recovering the quantum state proportional to $[\bb{U}, \bb{V}]^T$ is also of the same order as obtaining $[\tau^{-1}\bb{U}, \bb{V}]^T$.

To see why $\kappa(\bb{L}') \lesssim \tau^{-1}$ we use the following argument. Let $\bb{E} = \text{diag}(\tau^{-1}\bb{I}, \bb{I})$, where $\bb{I}$ is the identical matrix. Then  $\bb{L}' = \tilde{\bb{L}} \bb{E}$. According to the bounds established later in Theorem \ref{thm:diffusive}, we have for small $\varepsilon$ that
$\sigma_{\max}(\tilde{\bb{L}}) \lesssim 1$ and $\sigma_{\min}(\tilde{\bb{L}}) \gtrsim \tau$.
For the maximum singular value, one has
$\sigma_{\max}(\bb{L}') = \|\bb{L}'\| \le \|\tilde{\bb{L}}\| \|\bb{E}\| \le \tau^{-1} \sigma_{\max}(\tilde{\bb{L}}) \lesssim \tau^{-1}$.
For the minimum singular value, one has $\sigma_{\min}(\bb{L}') = 1/\sigma_{\max}(\bb{L}'^{-1})$ by definition.
Noting that
$\sigma_{\max}(\bb{L}'^{-1}) = \|\bb{L}'^{-1}\| \le \|\bb{E}^{-1}\|\|\tilde{\bb{L}}^{-1}\|  \le \tau \sigma_{\max}(\tilde{\bb{L}}^{-1})
= \tau/\sigma_{\min}(\tilde{\bb{L}}) \lesssim 1$,
we obtain $\sigma_{\min}(\bb{L}') \gtrsim 1$.
Hence, $\kappa(\bb{L}') \lesssim \tau^{-1}$ when $\varepsilon$ is small, as for $\tilde{\bb{L}}$.

We further remark that the formulation in Eq.~\eqref{AP2Pre} can raise the success probability of obtaining $\bb{U}$ by a factor of $1/\tau^2$ (since $1/\tau>1$) if we obtain the output $(\tau^{-1}\bb{U}, \bb{V})$ instead of $(\bb{U}, \bb{V})$. Similarly, we can instead boost the success probability of obtaining $\bb{V}$ by a factor of $1/\tau^2$ by replacing \eqref{AP2Pre2} with
\begin{equation}
\bb{L}' \begin{bmatrix}
\bb{I}  &  \\
 & \tau \bb{I}
\end{bmatrix} \begin{bmatrix}
\bb{U}  \\
\tau^{-1} \bb{V}
\end{bmatrix} = \begin{bmatrix}
\tau^{-1}\bb{L}_{11}  & \ \bb{L}_{12} \\
 \bb{L}_{12}  & \tau \bb{L}_{22}
\end{bmatrix}
\begin{bmatrix}
\bb{U}  \\
\tau^{-1}\bb{V}
\end{bmatrix}
= \begin{bmatrix}
\tau^{-1}\bb{F}_1  \\
\bb{F}_2
\end{bmatrix},
\end{equation}
where the new matrix still satisfies the condition number $\lesssim \tau^{-1}$ for small $\varepsilon$.

Since the condition numbers and the max-norms of our matrices $\tilde{\bb{L}}$, $\bb{L}'$ and in Eq.~\eqref{AP2Pre2} are comparable in size, the query complexity in obtaining either $[\bb{U}, \bb{V}]$, $[\tau^{-1}\bb{U}, \bb{V}]$ or $[\bb{U}, \tau^{-1}\bb{V}]$ are the same. We make the additional remark that the max-norms of these matrices can scale like $O(1/\tau)$, but the contribution to the total query complexity will still be independent of $\varepsilon$ since $\tau$ is independent of $\varepsilon$. Thus our total query complexity in obtaining these states will also be independent of $\varepsilon$, showing an additional advantage of using AP schemes.

\subsubsection{A penalized diffusive relaxation scheme}

The AP schemes given previously tend to explicit discretizations of the heat equation when $\varepsilon\to 0$,  thus suffer from the standard parabolic CFL restriction: $\tau \le Ch^2$. To overcome this problem, one can  employ the penalization technique proposed in \cite{Russo-2013}. As an example, we only consider the diffusive relaxation scheme.

By adding $\mu u_{xx}$ to the left and right sides of the first equation of \eqref{relaxation}, where $\mu = \mu(\varepsilon) \to 1$ as $\varepsilon \to 0$, one obtains
\begin{align*}
\begin{cases}
u_t + (v+ \mu u_x)_x = \mu u_{xx}, \\
v_t + u_x = - \frac{1}{\varepsilon}\Big( v + (1-\varepsilon) u_x \Big), \quad 0<\varepsilon \ll 1.
\end{cases}
\end{align*}
The way to remove the parabolic CFL restriction is to treat $u_{xx}$ implicitly since when $\varepsilon \to 0$ it will give an implicit discretization of the heat equation. According to the construction of the diffusive relaxation scheme, the scheme can be written as
\begin{enumerate}
  \item Relaxation step
  \begin{equation} \label{relaxationstep1}
\begin{cases}
\frac{u_j^{*} - u_j^n}{\tau}  = \mu \frac{u_{j-1}^*-2u_j^*+u_{j+1}^*}{h^2}, \quad j = 1,\cdots, N_x-1, \\
\frac{v_j^{*} - v_j^n}{\tau}  = - \frac{1}{\varepsilon}\Big( v_j^{*} + (1-\varepsilon)\frac{u_{j+1}^{*}-u_{j-1}^{*}}{2h} \Big),
\quad j = 1,\cdots, N_x-1.
\end{cases}
\end{equation}
 In the computation, we take $v_0^* = v_0^n$ and $v_{N_x}^* = v_{N_x}^n$.
  \item Convection step
  \begin{equation} \label{convectionstep1}
\begin{cases}
\frac{u_j^{n+1} - u_j^*}{\tau} + \frac{v_{j+1}^*-v_{j-1}^*}{2h} - \frac{h}{2}\frac{u_{j-1}^*-2u_j^*+u_{j+1}^*}{h^2}+ \mu \frac{u_{j-1}^*-2u_j^*+u_{j+1}^*}{h^2} = 0, \\
\frac{v_j^{n+1} - v_j^*}{\tau} + \frac{u_{j+1}^*-u_{j-1}^*}{2h} - \frac{h}{2}\frac{v_{j-1}^*-2v_j^*+v_{j+1}^*}{h^2} = 0.
\end{cases}
\end{equation}
\end{enumerate}

In the following we set $\mu = 1$. Let $\widetilde{\bb{B}} = \bb{I} -\tilde{\beta}\bb{L}_h$, where $\tilde{\beta} = \beta/h = \tau/h^2$. Then \eqref{relaxationstep1} and \eqref{convectionstep1} can be written in vector form as
\[
\begin{cases}
 \widetilde{\bb{B}}\bb{u}^* = \bb{u}^n + \tilde{\beta}\bb{b}^n, \\
 \gamma\bb{v}^* = \bb{v}^n - \nu \bb{A}\bb{u}^* - \nu \frac{\beta}{2} \widetilde{\bb{b}}^n
\end{cases}
\]
and
\[
\begin{cases}
\bb{u}^{n+1} = (\bb{B}+\widetilde{\bb{B}}-\bb{I}) \bb{u}^* - \bb{A}\bb{v}^* + \frac{\beta}{2}(\bb{b}^n - \widetilde{\bb{c}}^n)-\tilde{\beta}\bb{b}^n, \\
\bb{v}^{n+1} = \bb{B} \bb{v}^* - \bb{A}\bb{u}^* + \frac{\beta}{2}(\bb{c}^n - \widetilde{\bb{b}}^n),
\end{cases}
\]
respectively. Introducing the following notations
\begin{align}
\bb{A}_1 = \frac{1}{\gamma}\bb{A} = \frac{\beta}{2\gamma}\bb{M}_h, \quad
\bb{A}_2 = \frac{1}{\gamma}\bb{B} = \frac{1}{\gamma} (\bb{I} + \frac{\beta}{2} \bb{L}_h), \label{Aeq} \\
\bb{B}_1 = (\widetilde{\bb{B}} - \bb{I} + \bb{B} +  \frac{\nu}{\gamma}\bb{A}^2 ) \widetilde{\bb{B}}^{-1}, \quad
\bb{B}_2 = (\bb{A} + \frac{\nu}{\gamma}\bb{BA})\widetilde{\bb{B}}^{-1} , \label{Beq}
\end{align}
\[\widetilde{\bb{f}}^{n+1} = \tilde{\beta}\bb{B}_1\bb{b}^n + \frac{\beta}{2}(\bb{b}^n - \widetilde{\bb{c}}^n) + \frac{\nu\beta}{2\gamma}\bb{A}\widetilde{\bb{b}}^n - \tilde{\beta}\bb{b}^n, \]
\[\widetilde{\bb{g}}^{n+1} = -\tilde{\beta}\bb{B}_2\bb{b}^n - \frac{\beta}{2}( \widetilde{\bb{b}}^n - \bb{c}^n) - \frac{\nu\beta}{2\gamma}\bb{B}\widetilde{\bb{b}}^n,\]
and eliminating $\bb{u}^*$ and $\bb{v}^*$, it is easy to get the linear system
\begin{equation}\label{APPenalty}
\bb{L}_{\text{penalized}} \bb{S} = \bb{F}_{\text{penalized}},
\end{equation}
where $\bb{L}_{\text{penalized}} = (\bb{L}_{ij})_{2\times 2}$, $\bb{F}_{\text{penalized}} = [\bb{F}_1; \bb{F}_2]$, and
\[\bb{L}_{11} =
\begin{bmatrix}
\bb{I}  &            &           &            \\
-\bb{B}_1 & \bb{I}     &           &            \\
        &\ddots      & \ddots    &    \\
        &            & -\bb{B}_1   & \bb{I}     \\
\end{bmatrix}
, \qquad
\bb{L}_{12} =
\begin{bmatrix}
\bb{O}  &            &           &            \\
\bb{A}_1 & \bb{O}     &           &            \\
        &\ddots      & \ddots    &    \\
        &            & \bb{A}_1   & \bb{O}     \\
\end{bmatrix}, \qquad
\bb{F}_1 =
\begin{bmatrix}
\widetilde{\bb{f}}^1 + \bb{B}_1\bb{u}^0-\bb{A}_1\bb{v}^0 \\
\widetilde{\bb{f}}^2  \\
\vdots\\
\widetilde{\bb{f}}^{N_t} \\
\end{bmatrix}
\]
\[\bb{L}_{21} =
\begin{bmatrix}
\bb{O}   &            &            &            \\
\bb{B}_2 & \bb{O}     &            &            \\
         & \ddots     & \ddots     &    \\
         &            & \bb{B}_2   & \bb{O}     \\
\end{bmatrix}
, \qquad
\bb{L}_{22} =
\begin{bmatrix}
\bb{I}  &            &           &            \\
-\bb{A}_2 & \bb{I}     &           &            \\
        &\ddots      & \ddots    &    \\
        &            & -\bb{A}_2   & \bb{I}    \\
\end{bmatrix}
, \qquad
\bb{F}_2 =
\begin{bmatrix}
\widetilde{\bb{g}}^1 - \bb{B}_2\bb{u}^0 + \bb{A}_2\bb{v}^0\\
\widetilde{\bb{g}}^2  \\
\vdots\\
\widetilde{\bb{g}}^{N_t} \\
\end{bmatrix}.
\]

\subsection{Time complexity analysis of the AP schemes}

Denote $C$ and $Q$ to be the time complexity of classical and quantum processing, respectively. Fix the mesh size $h = 1/N_x$, and set $\delta$ to be the desired error bound of the quantum linear solver. Since the truncation error is $\mathcal{O}(\tau + h)$, we take $h = \mathcal{O}(\delta)$ in the following.

\subsubsection{IMEX scheme}

The previous qualitative analysis of the relationship between the matrix $\bb{H}$ and the scaling parameter $\varepsilon$ led to the design of the preconditioned IMEX scheme for solving the problem \eqref{heatstiff}. Now, we estimate the upper and lower bounds of the eigenvalues of $\bb{H}$, thus give a quantification of the time complexity.

To this end, we first list several lemmas.

\begin{lemma}\cite[Eq. (2.3)]{Rump-2011}\label{lem:IE}
Let $\bb{E} \in \mathbb{R}^{n\times n}$ and $\|\bb{E}\| \le \alpha<1$. Then $\bb{I} + \bb{E}$ is invertible and its singular values satisfy
$1 - \alpha \le \sigma( \bb{I} + \bb{E} ) \le 1 + \alpha$.
\end{lemma}

\begin{lemma}\cite[Weyl's inequality]{Stewart-1991}\label{lem:Weyl}
Let $ \bb{A}$ be a square matrix of order $n$ with singular values $\sigma_1\ge \sigma_2 \ge \cdots \sigma_n$. Denote $\tilde{\bb{A}} = \bb{A} + \bb{E}$ to be a perturbation of $\bb{A}$ with singular values $\tilde{\sigma}_1\ge \tilde{\sigma}_2 \ge \cdots \tilde{\sigma}_n$. Then there holds
$|\tilde{\sigma}_i - \sigma_i| \le \|\bb{E}\|$ for $i = 1,2,\cdots, n$.
\end{lemma}

\begin{theorem} \label{thm:IMEX}
Let $\beta = \tau/h$ and $\tau, h\le 1$. The eigenvalue of the matrix $\bb{H}$ corresponding to \eqref{AP1Pre} is denoted by $\lambda$.
\begin{enumerate}[(1)]
  \item For all $\varepsilon \ge 0$, $|\lambda|_{\max} \le 2 + \beta$. If taking $\tau = \iota h^2$ and $\beta = \iota h\le 1/2$ , where $\iota$ is a constant, then
  \[ |\lambda|_{\min} \ge \frac{h^2}{8} - (3\beta + 2-\tau) \frac{\varepsilon}{1+\varepsilon}.\]
  Hence for $\varepsilon \to 0$ the condition number and the sparsity of $\bb{H}$ satisfy
  \[\kappa = \mathcal{O}(N_t) \qquad \mbox{and} \qquad s = \mathcal{O}(1).\]
  \item  Let $\varepsilon \to 0$ and take $\tau = \mathcal{O}(h^2)$. Then the time complexity of the IMEX scheme is
  \[C_{\text{IMEX}} = \mathcal{O}(N_x^3) = \mathcal{O}(\delta^{-3}), \qquad Q_{\text{IMEX}} = \mathcal{O}(N_x^2 \log N_x)
  = \mathcal{O}(\delta^{-2} \log (1/\delta)).\]
\end{enumerate}
\end{theorem}
\begin{proof}
1) We first bound $|\lambda|_{\max}$. The coefficient matrix of \eqref{AP1Pre} is $\widetilde{\bb{L}}_{\text{IMEX}} = (l_{ij})$. For convenience, we omit the subscript. Denote the order by $n$. Let $\mu$ be the eigenvalue of $\widetilde{\bb{L}}\widetilde{\bb{L}}^T$ and define
\[R_i = \sum\limits_{j=1}^n |l_{ij}|, \qquad C_j = \sum\limits_{i=1}^n |l_{ij}|.\]
We have $\mu \le  \max_{ij} R_i C_j$ since $\|\widetilde{\bb{L}}\|^2 \le \|\widetilde{\bb{L}}\|_1 \|\widetilde{\bb{L}}\|_\infty$.
A direct calculation shows that the possible values of $R_i$ are
\begin{align*}
& a_1 = 1,  \quad a_2 = 2,  \quad a_3 = 2+\beta, \quad a_4 = \varrho\nu\beta/2 + \varrho\gamma, \\
& a_5 = \varrho\nu\beta + \varrho\gamma,  \quad
  a_6 = (\varrho\nu + \varrho)\beta+ \varrho + \varrho\gamma.
\end{align*}
By the definition of the notations,
\[\varrho = \frac{\varepsilon}{1+\varepsilon}<1, \quad
\varrho\nu + \varrho = \frac{1}{1+\varepsilon} <1, \quad
\varrho \gamma = \frac{\tau+\varepsilon}{1+\varepsilon}<1, \]
which implies
\[r_i \le a_3 = 2+\beta, \quad i \ge 1.\]
Similar calculation shows that the possible values of $C_j$ are
\begin{align*}
& b_1 = \varrho\gamma,  \quad b_2 = 1+\varrho\nu\beta/2,  \quad b_3 = 1+\varrho\nu\beta, \quad
  b_4 = \beta + \varrho\gamma + \varrho, \\
& b_5 = 2 + (\varrho\nu + \varrho)\beta,  \quad
  b_6 = \varrho + \varrho\gamma + (1-\varrho)\beta/2, \quad
  b_7 = 2 + (\varrho + \varrho\nu-1)\beta/2.
\end{align*}
Obviously,
\[C_j \le 2+\beta, \quad j \ge 1.\]
Combining the estimates of $R_i$ and $C_j$, we have $\mu \le (2+\beta)^2$, and hence $|\lambda|_{\max} \le 2 + \beta$.

2) We first briefly explain the idea of bounding $|\lambda|_{\min}$: Let $\bb{L}_{\varepsilon}$ be the coefficient matrix and $\bb{L}_{\varepsilon} = \bb{L}_0 +\bb{E}$, where $\bb{L}_0$ is the coefficient matrix with $\varepsilon=0$. By the Weyl's inequality in Lemma \ref{lem:Weyl}, it suffices to determine the lower bound of $|\lambda(\bb{L}_0)|_{\min}$ and the upper bound of $\|\bb{E}\|$.

For this reason, we first let $\varepsilon=0$ and for convenience assume that $\widetilde{\bb{L}}_{ij}$ has only three row or column blocks. Let $\widetilde{\bb{L}} = \bb{I} + \bb{E}$, where
\[ \bb{E}  = {\scriptsize \left[\begin{array}{ccc:cccccccc}
\bb{O}   & \bb{O}   & \bb{O}    & \bb{O}  & \bb{O}  & \bb{O} \\
-\bb{B}  & \bb{O}   & \bb{O}    & \bb{A}  & \bb{O}  & \bb{O} \\
 \bb{O}  & -\bb{B}  & \bb{O}    & \bb{O}  & \bb{A}  & \bb{O} \\
 \hdashline
\bb{A}   & \bb{O}   & \bb{O}    & -d\bb{I}  & \bb{O}  & \bb{O} \\
\bb{O}   & \bb{A}   & \bb{O}    & \bb{O}  & -d\bb{I}  & \bb{O} \\
 \bb{O}  &  \bb{O}  & \bb{A}    & \bb{O}  & \bb{O}  & -d\bb{I} \\
\end{array} \right]}, \qquad d = 1 - \tau.\]
We now estimate $\alpha:=\|\bb{E}\| = \lambda_{\max}(\bb{E}\bb{E}^T)^{1/2}$. A direct manipulation gives
\[\scriptsize \bb{E}\bb{E}^T = \left[\begin{array}{ccc:cccccccc}
\bb{O}   &                                 &      & \bb{O}                   &         &   \\
         & \bb{B}\bb{B}^T+\bb{A}\bb{A}^T   &      & -\bb{B}\bb{A}^T-d\bb{A}  & \bb{O}  &   \\
         &          & \bb{B}\bb{B}^T+\bb{A}\bb{A}^T    & \bb{O}  & -\bb{B}\bb{A}^T-d\bb{A}  & \bb{O} \\
 \hdashline
   \bb{O}      & -\bb{A}\bb{B}^T-d\bb{A}^T   &     & \bb{A}\bb{A}^T+d^2\bb{I}  &    &   \\
         &   \bb{O}      & -\bb{A}\bb{B}^T-d\bb{A}^T    &    & \bb{A}\bb{A}^T+d^2\bb{I}  &   \\
         &         &    \bb{O}    &    &    & \bb{A}\bb{A}^T+d^2\bb{I} \\
\end{array} \right].
\]
It is apparent that the eigenvalues of $\bb{A}\bb{A}^T+d^2\bb{I}$ are the ones of $\bb{E}\bb{E}^T$, given as
\[d^2 + \lambda(\bb{A}\bb{A}^T) = d^2 - \lambda(\bb{A}^2). \]
Notice that the eigenvalues of $\bb{M}_h = {\rm tril}(-1, 0, 1)$ are zero or pure imaginary, thus $-\lambda(\bb{M}_h^2)\ge 0$. One can check that $8h^2 < |\lambda(\bb{M}_h^2)| <4$. Since $\bb{A} = \frac{\beta}{2} \bb{M}_h$ and $d = 1-\tau$,
\[(1-\tau)^2 + 2\beta^2 h^2 < d^2 + \lambda(\bb{A}\bb{A}^T) < (1-\tau)^2 + \beta^2.\]

It remains to consider the middle part of $\bb{E}\bb{E}^T$, i.e.,
\[\scriptsize \widetilde{\bb{E}}\widetilde{\bb{E}}^T = \left[\begin{array}{cc:cccc}
 \bb{B}\bb{B}^T+\bb{A}\bb{A}^T   &      & -\bb{B}\bb{A}^T-d\bb{A}  &   \\
                & \bb{B}\bb{B}^T+\bb{A}\bb{A}^T  &    & -\bb{B}\bb{A}^T-d\bb{A}  \\
 \hdashline
 -\bb{A}\bb{B}^T-d\bb{A}^T   &     & \bb{A}\bb{A}^T+d^2\bb{I}  &      \\
        & -\bb{A}\bb{B}^T-d\bb{A}^T    &    & \bb{A}\bb{A}^T+d^2\bb{I}    \\
\end{array} \right],
\]
where,
\[\scriptsize \widetilde{\bb{E}} = \left[\begin{array}{cc:cccc}
 -\bb{B}   &         &     \bb{A}   &   \\
           & -\bb{B} &              & \bb{A}  \\
 \hdashline
  \bb{A}   &          & -d\bb{I}  &      \\
           & \bb{A}   &           & -d\bb{I}    \\
\end{array} \right].
\]
Let
\[\scriptsize \widehat{\bb{E}} = \left[\begin{array}{cc:cccc}
 -\bb{I}   &         &        &   \\
           & -\bb{I} &              &  \\
 \hdashline
             &          & \bb{I}  &      \\
           &    &           & \bb{I}    \\
\end{array} \right] \left[\begin{array}{cc:cccc}
 -\bb{B}   &         &     \bb{A}   &   \\
           & -\bb{B} &              & \bb{A}  \\
 \hdashline
  \bb{A}   &          & -d\bb{I}  &      \\
           & \bb{A}   &           & -d\bb{I}    \\
\end{array} \right] = \left[\begin{array}{cc:cccc}
 \bb{B}   &         &     -\bb{A}   &   \\
           & \bb{B} &              & -\bb{A}  \\
 \hdashline
  \bb{A}   &          & -d\bb{I}  &      \\
           & \bb{A}   &           & -d\bb{I}    \\
\end{array} \right].
\]
Then $\widehat{\bb{E}}$ is symmetric and $\|\widetilde{\bb{E}}\| = \|\widehat{\bb{E}}\|$ since the orthogonal transformations do not change the singular values. The matrix $\widehat{\bb{E}}$ can be decomposed as $\widehat{\bb{E}} = \bb{I} + \frac{\beta}{2} \bb{C}_h$ with
\[\bb{C}_h = \left[\begin{array}{cc:cccc}
 \bb{L}_h   &         &     -\bb{M}_h   &   \\
           & \bb{L}_h &              & -\bb{M}_h  \\
 \hdashline
 \bb{M}_h   &          & -\hat{d}\bb{I}  &      \\
           & \bb{M}_h   &           & -\hat{d}\bb{I}    \\
\end{array} \right], \qquad \hat{d} = \frac{2(2-\tau)}{\beta}\]
consisting of two basic matrices $\bb{L}_h$ and $\bb{M}_h$ given in \eqref{MhLh}. The matrix $\bb{C}_h$ can be diagonalized and has real eigenvalues. In fact, $\bb{C}_h$ is negative definite. By definition, for any $\bb{x} \ne \bb{0}$ we write it as $\bb{x} = [\bb{x}_1; \bb{x}_2; \bb{x}_3; \bb{x}_4]$. The Cauchy-Schwarz inequality yields
\begin{align*}
\bb{x}^T \bb{C}_h \bb{x}
& = \bb{x}_1^T\bb{L}_h \bb{x}_1 + \bb{x}_2^T \bb{L}_h\bb{x}_2
   + 2\bb{x}_3^T \bb{M}_h\bb{x}_1 + 2\bb{x}_4^T\bb{M}_h\bb{x}_2 - \hat{d}\bb{x}_3^T\bb{x}_3 - \hat{d}\bb{x}_4^T\bb{x}_4 \\
& \le \bb{x}_1^T(\bb{I}+\bb{L}_h) \bb{x}_1 + \bb{x}_2^T (\bb{I}+\bb{L}_h)\bb{x}_2
   -\bb{x}_3^T\bb{M}_h^2\bb{x}_3 - \bb{x}_4^T\bb{M}_h^2\bb{x}_4 - \hat{d}\bb{x}_3^T\bb{x}_3 - \hat{d}\bb{x}_4^T\bb{x}_4 \\
& \le \Big(1 + \lambda_{\max}(\bb{L}_h) +\lambda_{\max}(- \bb{M}_h^2) - \hat{d}\Big)\bb{x}^T\bb{x}.
\end{align*}
Obviously,
\[1 + \lambda_{\max}(\bb{L}_h) \bb{x}^T\bb{x} +\lambda_{\max}(- \bb{M}_h^2) - \hat{d} < 5-\hat{d} - 8h^2
 = - \frac{4}{\iota h} + 2h - 8h^2 + 5,\]
which shows that when $\beta = \iota h \le 1/2$ and $h\le 1$, the eigenvalues of $\bb{C}_h$ are less than zero. Hence,
\[\lambda(\widehat{\bb{E}}) = \lambda(\bb{I} + \frac{\beta}{2}\bb{C}_h) < -1 + \frac{5}{2}\iota h + \iota h^2 - 4\iota h^3.\]
By the Rayleigh quotient theorem \cite{Horn2013} and noting that $(\bb{L}_h)_{ii} = -2 ~(i\ge 1)$,
\[\lambda_{\min} (\bb{C}_h) \le \min_{i} (\bb{C}_h)_{ii} = \min \{-2,  -\hat{d} \} = - \hat{d},\]
and hence $\lambda(\widehat{\bb{E}}) \ge  -1+\iota h^2$.
Combining the above equations, we then obtain $|\lambda(\widehat{\bb{E}})| \le  1 - \iota h^2$.

A simple calculation shows that for $h<1$ it holds $1 - \iota h^2 \le \sqrt{(1-\tau)^2 + \beta^2 }$, where $\beta = \iota h \le 1/2$. Summing up the previous analysis, we get $\alpha \le \sqrt{(1-\tau)^2 + \beta^2 }$, which in turn gives
\begin{align*}
1 - \alpha
& > \frac{1 -(1-\tau)^2 - \beta^2 }{1 +\sqrt{(1-\tau)^2 + \beta^2 } } = \frac{(1-2\iota)h^2 - \iota^2h^4}{1 +\sqrt{(1-\tau)^2 + \beta^2 } } \\
& = \frac{(1-\iota)^2h^2 + \iota^2h^2(1-h^2)}{1 +\sqrt{(1-\tau)^2 + \beta^2 } } \ge \frac{h^2}{8}.
\end{align*}
The desired estimate is obtained by using Lemma \ref{lem:IE}.

3) For $\varepsilon>0$, let $\widetilde{\bb{L}}^0$ be the matrix corresponding to $\varepsilon=0$ and denote by $\tilde{\sigma}_i$ the singular values. Let $\bb{E} = \widetilde{\bb{L}} - \widetilde{\bb{L}}^0 = (\bb{E}_{ij})_{2\times 2}$, where $\bb{E}_{11} = \bb{E}_{12} = \bb{O}$, and
\[\bb{E}_{21} =
\begin{bmatrix}
(\varrho\nu-1)\bb{A}  &            &           &            \\
\varrho\bb{A} & (\varrho\nu-1)\bb{A}     &           &            \\
        &\ddots      & \ddots    &    \\
        &            & \varrho\bb{A}   & (\varrho\nu-1)\bb{A}     \\
\end{bmatrix},\]
\[
\bb{E}_{22} =
\begin{bmatrix}
(\varrho\gamma - \tau)\bb{I}  &            &           &            \\
-\varrho\bb{B} & (\varrho\gamma - \tau)\bb{I}     &           &            \\
        &\ddots      & \ddots    &    \\
        &            & -\varrho\bb{B}   & (\varrho\gamma - \tau)\bb{I}    \\
\end{bmatrix}.
\]
By Lemma \ref{lem:Weyl}, $|\sigma_i - \tilde{\sigma}_i| \le \|\bb{E}\|$,
where $\sigma_i$ are the singular values of $\bb{\tilde{L}}$. Let $n=N_x-1$. A direct calculation gives
\[
R_i(\bb{A}) = C_i(\bb{A}) = \begin{cases}
\frac{1}{2}\beta, \quad & i = 1,n\\
\beta, \quad & i = 2, \cdots, n-1
\end{cases}, \quad
R_i(\bb{B}) = C_i(\bb{B}) = \begin{cases}
1-\frac{1}{2}\beta, \quad & i = 1,n\\
1, \quad & i = 2, \cdots, n-1
\end{cases}.
\]
It follows from Lemma \ref{lem:Qi} that
\begin{align*}
\|\bb{E}\|
& \le ( \varrho + |1-\varrho\nu| ) R_{\max}(\bb{A} ) + ( \varrho R_{\max}(\bb{B}) + |\varrho\gamma - \tau| ) \\
& = \Big( \frac{\delta}{1+\delta} +  \frac{2\delta}{1+\delta} \Big) \beta +
\Big( \frac{\delta}{1+\delta}  + \frac{\delta}{1+\delta}(1-\tau)  \Big) \\
& = (3\beta + 2-\tau) \frac{\varepsilon}{1+\varepsilon},
\end{align*}
which naturally leads to the estimate for $|\lambda|_{\min}$.

4) Now we analyze the time complexity. The classical treatment is to iteratively solve \eqref{imex0}, i.e.,
\[
\begin{cases}
-\bb{B}\bb{u}^n + \bb{u}^{n+1} + \bb{A}\bb{v}^n = \bb{f}^{n+1}, \\
-\bb{B}\bb{v}^n + \gamma \bb{v}^{n+1} + \bb{A}\bb{u}^n + \nu \bb{A}\bb{u}^{n+1} = \bb{g}^{n+1}.
\end{cases}
\]
Given $\bb{u}^n$, $\bb{v}^n$ and $\bb{f}^{n+1}$, let $n_u=N_x-1$. Then the number of fundamental operations in obtaining $\bb{u}^{n+1}$ is
$\mathcal{O}(n_us_A + n_us_B) = \mathcal{O}(N_x)$.
With $\bb{u}^{n+1}$, we can compute $\bb{v}^{n+1}$ explicitly from the second equation, and the number of basic operations is
$\mathcal{O}(2n_us_A + ns_B) = \mathcal{O}(N_x)$.
That is, the time complexity of each iteration step is $\mathcal{O}(N_x)$, and thus the time complexity after $N_t$ iterations is
\[C_{\text{IMEX}} = \mathcal{O}(N_t N_x) = \mathcal{O}(N_x^3),\]
where $\tau = \mathcal{O}(h^2)$ is used in the last step.

For the quantum treatment, by the estimates of eigenvalues, the condition number $\kappa = \mathcal{O}(h^{-2}) = \mathcal{O}(N_x^2)$. Plugging in \eqref{cpCAS} yields
\[Q_{\text{IMEX}} = \mathcal{O}(N_x^2 \log(1/\delta)) = \mathcal{O}(N_x^2 \log N_x),\]
as required.
\end{proof}

\begin{figure}[!ht]
  \centering
  \includegraphics[scale=0.6]{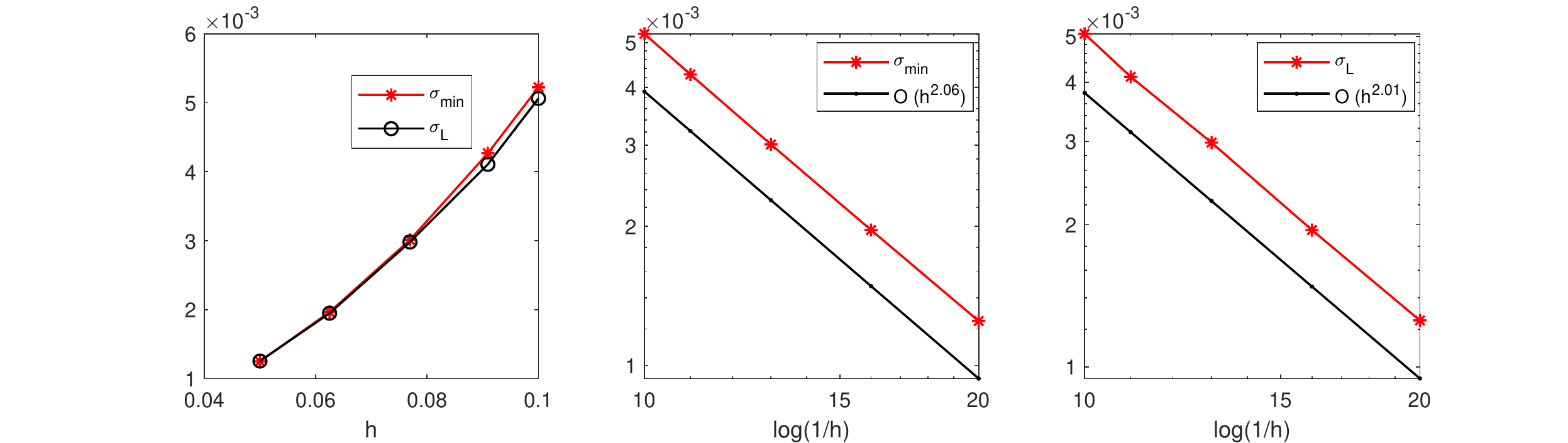}\\
  \caption{The estimated value $\sigma_L$ and the true value $\sigma_{\min}$ of the minimum singular value for the IMEX scheme}\label{fig:singularIMEX}
\end{figure}

The estimate of the minimum eigenvalue is optimal, and the best choice of $\alpha$ in the proof is $\alpha = \sqrt{(1 - \tau)^2 + \beta^2 }$. Let the mesh size be chosen uniformly in $[0.05, 0.1]$ with five points. The time step is taken as $\tau = h^2$, and the scaling parameter is $\varepsilon = 10^{-8}$. In Fig.~\ref{fig:singularIMEX}, we display the estimated and true values of the minimum eigenvalue, where $\sigma_L = 1-\alpha$ is the estimated value, which is close to the true value $\sigma_{\min}$ and always smaller than the true value. It is obvious that both values are second order with respect to $h$.

\subsubsection{Diffusive relaxation scheme}

We now consider the diffusive relaxation scheme \eqref{APsplit}.

\begin{theorem}\label{thm:diffusive}
Let $\beta = \tau/h\le 1$ and $\tau,h\le 1$. The eigenvalue of $\bb{H}$ for the reformulated diffusive relaxation scheme \eqref{AP2Pre} is denoted by $\lambda$.
\begin{enumerate}[(1)]
  \item If $\tau/h^2 \le 4$, then for $\varepsilon\ge 0$ there holds
  \[|\lambda|_{\max} \le 5 + \alpha, \qquad |\lambda|_{\min} \ge \frac{\tau}{4} - \alpha,\]
where
\[\alpha = \frac{\varepsilon}{\tau+\varepsilon}\Big( \frac{1+\tau}{\tau}\frac{\beta^2}{2} + \frac{\beta}{\tau} +
  (1+\tau)(\beta-\frac{\beta^2}{2}) + \tau \Big).\]
  Hence for $\varepsilon \to 0$ the condition number and the sparsity of $\bb{H}$ satisfy
  \[\kappa = \mathcal{O}(N_t) \qquad \mbox{and} \qquad s = \mathcal{O}(1).\]
    \item If taking $\tau = \mathcal{O}(h^2)$, then the time complexity of the diffusive relaxation scheme is
  \[C_{\text{relaxation}} = \mathcal{O}(N_x^3) = \mathcal{O}(\delta^{-3}), \qquad Q_{\text{relaxation}} = \mathcal{O}(N_x^2 \log N_x)
  = \mathcal{O}(\delta^{-2} \log(1/\delta)).\]
\end{enumerate}
\end{theorem}

\begin{proof}
(1) For simplicity, we omit the subscript ``relaxation'' and denote $\bb{L}_{\varepsilon}$ and $\bb{L}$ to be the reformulated coefficient matrix in \eqref{AP2Pre} with $\varepsilon\ge 0$ and $\varepsilon=0$, respectively. Let $\bb{L}_{\varepsilon} = \bb{L} + \bb{E}$. By the Weyl's inequality in Lemma \ref{lem:Weyl}, it suffices to bound $\sigma(\bb{L})$ and $\|\bb{E}\|$, where
\[ \bb{L}  = \begin{bmatrix}
\bb{L}_{11} & \bb{O} \\
\bb{L}_{21} & \bb{I}
\end{bmatrix} = {\scriptsize \left[\begin{array}{cccc:cccccccc}
\bb{I}    &          &              &         &   \bb{O}   &        &              &  \\
-\bb{B}_1 & \bb{I}   &              &         &            & \bb{O} &              &   \\
          & \ddots   & \ddots       &         &            &        &    \ddots    & \\
          &          & -\bb{B}_1    & \bb{I}  &            &        &              & \bb{O}\\
 \hdashline
\bb{O}     &         &          &        & \bb{I}    &          &          & \\
\bb{B}_2   & \bb{O}  &          &        &           & \bb{I}   &          & \\
           & \ddots  & \ddots   &        &           &          & \ddots   & \\
           &         & \bb{B}_2 & \bb{O} &           &          &          &  \bb{I} \\
\end{array} \right]},\]
with $\bb{B}_1 = \bb{B}+\tau^{-1}\bb{A}^2$ and $\bb{B}_2 = \tau \bb{A} + \bb{BA}$ since $\nu/\gamma\to \tau^{-1}$ as $\varepsilon\to 0$.

Step 1: When $\beta\le 1$, a direct calculation gives
\[
R_i(\bb{B}_2) = C_i(\bb{B}_2) = \begin{cases}
( \tau + (1-\beta) )\frac{\beta}{2} + \frac{\beta^2}{2}, \quad & i = 1,n,\\
( \tau + (1-\beta) )\beta + \frac{\beta^2}{4}, \quad & i = 2, n-1, \\
( \tau + (1-\beta) )\beta + \frac{\beta^2}{2}, \quad & i = 3,\cdots, n-2,
\end{cases}
\]
hence
\begin{equation}\label{rcB2}
\|\bb{B}_2\| \le  ( \tau + (1-\beta) )\beta + \frac{\beta^2}{2} \le 3.
\end{equation}
When $\tau/h^2 \le 4$, one has $0\le 1-\frac{\beta^2}{4\tau} = 1 - \frac{\tau}{4h^2} \le 1$. A similar calculation gives
\begin{equation}\label{RB1}
R_i(\bb{B}_1) = C_i(\bb{B}_1) = \begin{cases}
1-\frac{1}{2}\beta, \quad & i = 1, n,\\
1-\frac{\beta^2}{4\tau}, \quad & i = 2, n-1, \\
1, \quad & i = 3,\cdots, n-2,
\end{cases}
\end{equation}
and hence
\[\|\bb{B}_1\| \le  1, \qquad \|\bb{L}\| \le 5.\]

Step 2: It is obvious that
\[\begin{bmatrix}
\bb{L}_{11} & \bb{O} \\
\bb{L}_{21} & \bb{I}
\end{bmatrix}
\begin{bmatrix}
\bb{I} &  \bb{O} \\
-\bb{L}_{21} & \bb{I}
\end{bmatrix}=
\begin{bmatrix}
\bb{L}_{11} &  \\
            & \bb{I}
\end{bmatrix} \quad \text{or simply written as} \quad \bb{L}\bb{P} = \widehat{\bb{L}}.
\]
For any two matrices $\bb{A}$ and $\bb{B}$, there holds $\sigma_{\min}(\bb{AB}) \le \sigma_{\min}(\bb{A})\|\bb{B}\|$ (see Eq. (2.8) in \cite{Rump-2011}), which together with the previous equation gives
\begin{equation}\label{sigmge}
\sigma_{\min}(\bb{L}) \ge \sigma_{\min}(\widehat{\bb{L}})/\|\bb{P}\| \ge \sigma_{\min}(\bb{L}_{11})/\|\bb{P}\|.
\end{equation}
One easily finds that
\[\|\bb{P}\| \le 1 + \|\bb{B}_2\| < 4,\]
yielding
\begin{equation}\label{sigL}
\sigma_{\min}(\bb{L}) \ge  \sigma_{\min}(\bb{L}_{11})/4.
\end{equation}
Therefore it suffices to bound $\sigma_{\min}(\bb{L}_{11})$.

Step 3: Let $\bb{D} = \bb{L}_{11}\bb{L}_{11}^T$. Noting that $\bb{B}_1$ is symmetric, one has
\[\bb{D} =
\begin{bmatrix}
\bb{I}    &  -\bb{B}_1           &                      &            \\
-\bb{B}_1 & \bb{I} + \bb{B}_1^2  &  \ddots          &            \\
          &  \ddots              &  \ddots               & -\bb{B}_1    \\
          &                      &  -\bb{B}_1       &  \bb{I} + \bb{B}_1^2     \\
\end{bmatrix}.
 \]
Let $\bb{B}_1$ be similar to the diagonal matrix $\Lambda$ consisting of the eigenvalues. Then,
\[\bb{D} \sim \begin{bmatrix}
\bb{I}    &  -\Lambda           &                      &            \\
-\Lambda & \bb{I} + \Lambda^2  &  \ddots          &            \\
          &  \ddots              &  \ddots               & -\Lambda    \\
          &                      &  -\Lambda       &  \bb{I} + \Lambda^2     \\
\end{bmatrix} =: \widehat{\bb{D}}.
\]
We point out that the Gershgorin circle theorem only gives a trivial lower bound 0.
The matrix $\widehat{\bb{D}}$ has the same structure of $\bb{D}$.
Hence it can be written as $\widehat{\bb{D}} = \bb{C}\bb{C}^T$, where
\[\bb{C} =
\begin{bmatrix}
\bb{I}         &                 &            &           &       \\
-\Lambda       & \bb{I}          &            &           &        \\
               &  \ddots         &  \ddots    &           &     \\
               &                 &  \ddots    &  \ddots   &      \\
               &                 &            &  -\Lambda &  \bb{I}    \\
\end{bmatrix}, \qquad
\bb{C}^{-1} =
\begin{bmatrix}
\bb{I}          &              &                 &          &            \\
 \Lambda        & \bb{I}       &                 &          &        \\
 \Lambda^2      &  \ddots      &  \ddots         &          &\\
   \vdots       &  \ddots      &  \ddots         &  \ddots  &   \\
\Lambda^{N_t-1} &   \cdots     &   \Lambda^2     &  \Lambda &  \bb{I}      \\
\end{bmatrix}.
 \]
By the definition of the 2-norm,
\[\frac{1}{\lambda_{\min}(\bb{D})} = \frac{1}{\lambda_{\min}(\widehat{\bb{D}})}
= \|\widehat{\bb{D}}^{-1}\| = \|(\bb{C}\bb{C}^T)^{-1}\| \le \|\bb{C}^{-1}\|^2 = \sigma_{\max}^2(\bb{C}^{-1}).\]
Let $d$ be the maximum eigenvalue of $\bb{B}_1$. Since $\bb{B}_1$ is symmetric, $d=\|\bb{B}_1\| \le 1$. We obtain from Lemma \ref{lem:Qi} and the structure of $\bb{C}^{-1}$ that
\[\sigma_{\max}(\bb{C}^{-1}) \le 1 + d + d^2 + \cdots + d^{N_t-1} \le N_t = 1/\tau,\]
which combining \eqref{sigL} gives
\[\sigma_{\min}(\bb{L}) \ge \tau/4.\]

Step 4: We now estimate $\|\bb{E}\|$, where
\begin{equation}\label{EAP2}
 \bb{E} = {\scriptsize \left[\begin{array}{cccc:cccccccc}
\bb{O}    &          &               &         &   \bb{O}   &        &              &  \\
(\nu/\gamma-\tau^{-1})\bb{A}^2 & \bb{O}   &          &         &      \tau^{-1}/\gamma\bb{A}      & \bb{O} &              &   \\
          & \ddots   & \ddots       &         &            &        \ddots &    \ddots    & \\
          &          & (\nu/\gamma-\tau^{-1})\bb{A}^2    & \bb{O}  &            &        &         \tau^{-1}/\gamma\bb{A}     & \bb{O}\\
 \hdashline
\bb{O}     &         &          &        & \bb{O}    &          &          & \\
(\tau\nu/\gamma-1)\bb{BA}   & \bb{O}  &          &        &      -\tau/\gamma\bb{B}     & \bb{O}   &          & \\
           & \ddots  & \ddots   &        &           &      \ddots    & \ddots   & \\
           &         & (\tau\nu/\gamma-1)\bb{BA} & \bb{O} &           &          &      -\tau/\gamma\bb{B}    &  \bb{O} \\
\end{array} \right]},
\end{equation}
with
\[\frac{1}{\tau}-\frac{\nu}{\gamma} = \frac{\varepsilon}{\tau+\varepsilon}\frac{1+\tau}{\tau}, \quad
\frac{1}{\gamma} = \frac{\varepsilon}{\tau+\varepsilon}, \quad
1-\frac{\tau\nu}{\gamma}= \frac{\varepsilon}{\tau+\varepsilon}(1+\tau), \quad
\frac{\tau}{\gamma}= \frac{\varepsilon}{\tau+\varepsilon}\tau.\]
When $\beta\le 1$, one easily obtains
\[R_{\max}(\bb{A}^2) = \frac{\beta^2}{2}, \quad R_{\max}(\bb{A}) = \beta, \quad
R_{\max}(\bb{BA}) = \beta - \frac{\beta^2}{2}, \quad R_{\max}(\bb{B}) = 1. \]
Applying the Gershgorin circle theorem to get
\begin{align*}
\|\bb{E}\|
& \le \Big(\frac{1}{\tau}-\frac{\nu}{\gamma}\Big)R_{\max}(\bb{A}^2) + \frac{\tau^{-1}}{\gamma}R_{\max}(\bb{A})
+ \Big(1-\frac{\tau\nu}{\gamma}\Big) R_{\max}(\bb{BA}) + \frac{\tau}{\gamma} R_{\max}(\bb{B}) \\
& \le \frac{\varepsilon}{\tau+\varepsilon}\Big( \frac{1+\tau}{\tau}\frac{\beta^2}{2} + \frac{\beta}{\tau} +
  (1+\tau)(\beta-\frac{\beta^2}{2}) + \tau \Big).
\end{align*}
The desired estimate follows by applying the Weyl's inequality in Lemma \ref{lem:Weyl}.

(2) The time complexity of the diffusive relaxation scheme can be analyzed in the same manner as in Theorem \ref{thm:IMEX}.
\end{proof}

\subsubsection{The penalized diffusive relaxation scheme}

\begin{theorem}\label{thm:penalty}
Let $\beta = \tau/h \le 1/2$. Denote by $\lambda$ the eigenvalue of $\bb{H}$ for the penalized diffusive relaxation scheme \eqref{APPenalty}.
\begin{enumerate}[(1)]
  \item For every $\varepsilon\ge 0$
  \[|\lambda|_{\max} \le 4 + 2\alpha, \qquad  |\lambda|_{\min} \ge  \frac{\tau}{3} - \alpha,\]
where
\[\alpha = \frac{\varepsilon}{\tau+\varepsilon} \frac{2+4\tau}{\tau}.\]
Hence for $\varepsilon \to 0$ the condition number and the sparsity of $\bb{H}$ satisfy
  \[\kappa = \mathcal{O}(N_t) \qquad \mbox{and} \qquad s = \mathcal{O}(1).\]
  \item Let $\varepsilon \to 0$. If taking $\tau = \mathcal{O}(h)$ and assuming that $h$ is sufficiently small, then the time complexity of the penalized diffusive relaxation scheme is
  \[C_{\text{penalized}} = \mathcal{O}(N_x^{2.5}\log N_x) = \mathcal{O}(\delta^{-2.5}\log (1/\delta)), \quad
  Q_{\text{penalized}} = \mathcal{O}(N_x \log N_x) = \mathcal{O}(\delta^{-1} \log (1/\delta)).\]
  If $\tau = \mathcal{O}(h^2)$, then
  \[C_{\text{penalized}} = \mathcal{O}(N_x^3 \log N_x ) = \mathcal{O}(\delta^{-3}\log (1/\delta)), \quad Q_{\text{penalized}} = \mathcal{O}(N_x^2 \log N_x ) = \mathcal{O}(\delta^{-2}\log (1/\delta)).\]
\end{enumerate}
\end{theorem}

\begin{proof}
(1) Since the problem is linear, we can simply apply the discrete Fourier transform to characterize the singular values of the coefficient matrix $\bb{L}_{\text{penalized}}$ in \eqref{APPenalty}. For simplicity, we omit the subscript.

Step 1: To this end, we introduce the following expressions:
\begin{equation}\label{fouriervariables}
u_j^n = \hat{u}^n{\rm e}^{{\rm i} jkh}, \quad v_j^n = \hat{v}^n{\rm e}^{{\rm i} jkh},
\quad u_j^* = \hat{u}^*{\rm e}^{{\rm i} jkh}, \quad v_j^* = \hat{v}^*{\rm e}^{{\rm i} jkh},
\end{equation}
where $k$ represents the frequency variable and ${\rm i} = \sqrt{-1}$. Plugging them in \eqref{relaxationstep1}-\eqref{convectionstep1}, one obtains
\begin{enumerate}
  \item[] Relaxation step:
\[
\begin{cases}
\frac{\hat{u}^{*} - \hat{u}^n}{\tau}
=  \frac{({\rm e}^{-{\rm i} kh}-2 + {\rm e}^{{\rm i}kh})\hat{u}^*}{h^2}, \\
\frac{\hat{v}^{*} - \hat{v}^n}{\tau}
= - \frac{1}{\varepsilon}\Big( \hat{v}^{*} + (1-\varepsilon)\frac{({\rm e}^{{\rm i} kh}- {\rm e}^{-{\rm i}kh})\hat{u}^*}{2h} \Big).
\end{cases}
\]
  \item[] Convection step:
\[
\begin{cases}
\frac{\hat{u}^{n+1} - \hat{u}^*}{\tau} + \frac{({\rm e}^{{\rm i} kh}- {\rm e}^{-{\rm i}kh})\hat{v}^*}{2h} - \frac{h}{2}\frac{({\rm e}^{-{\rm i} kh}-2 + {\rm e}^{{\rm i}kh})\hat{u}^*}{h^2}
+ \frac{({\rm e}^{-{\rm i} kh}-2 + {\rm e}^{{\rm i}kh})\hat{u}^*}{h^2} = 0, \\
\frac{\hat{v}^{n+1} - \hat{v}^*}{\tau} + \frac{({\rm e}^{{\rm i} kh}- {\rm e}^{-{\rm i}kh})\hat{u}^*}{2h} - \frac{h}{2}\frac{({\rm e}^{-{\rm i} kh}-2 + {\rm e}^{{\rm i}kh})\hat{v}^*}{h^2} = 0.
\end{cases}
\]
\end{enumerate}
Collecting the above system, one has
\[\begin{cases}
\hat{u}^{n+1} + c_{1,\varepsilon} \hat{u}^n + c_{2,\varepsilon} \hat{v}^n = 0\\
\hat{v}^{n+1} + d_{1,\varepsilon} \hat{v}^n + d_{2,\varepsilon} \hat{u}^n = 0
\end{cases}, \quad n = 0,1,\cdots, N_t-1,\]
where
\begin{align*}
& c_{1,\varepsilon} = \frac{1}{1+4\tilde{\beta} \sin^2 \frac{kh}{2} } \Big( -1-4\tilde{\beta} \sin^2 \frac{kh}{2} + ( 2\beta + \frac{\beta^2\nu}{\gamma} ) \sin^2 (kh) \Big),\\
& c_{2,\varepsilon} = {\rm i} \frac{\beta}{\gamma} \sin (kh),  \qquad
  d_{1,\varepsilon} = \frac{2\beta}{\gamma} \sin^2 \frac{kh}{2}, \\
& d_{2,\varepsilon} = {\rm i} \frac{1}{1+4\tilde{\beta} \sin^2 \frac{kh}{2} } \Big( \beta  - \frac{2\beta^2\nu}{\gamma} \sin^2 (kh) \Big)\sin (kh),
\end{align*}
and
\[\gamma = 1+\tau/\varepsilon,  \qquad \nu = (1-\varepsilon)/\varepsilon.\]
The resulting coefficient matrix corresponding to $\bb{L}$ can be written as
\begin{equation}\label{LFourier}
\hat{\bb{L}}_{\varepsilon}  = \begin{bmatrix}
\hat{\bb{L}}_{11,\varepsilon} & \hat{\bb{L}}_{12,\varepsilon} \\
\hat{\bb{L}}_{21,\varepsilon} & \hat{\bb{L}}_{22,\varepsilon}
\end{bmatrix} = {\scriptsize \left[\begin{array}{cccc:cccccccc}
1         &          &              &         &   0      &            &              &  \\
c_{1,\varepsilon}       & 1        &              &         &      c_{2,\varepsilon}      &     0    &              &   \\
          & \ddots   & \ddots       &         &            &      \ddots       &    \ddots    & \\
          &          & c_{1,\varepsilon}          & 1       &            &            &      c_{2,\varepsilon}        & 0\\
 \hdashline
    0     &          &              &         & 1          &            &              & \\
d_{2,\varepsilon}       & 0        &              &         &  d_{1,\varepsilon}       & 1          &              & \\
          & \ddots   & \ddots       &         &            &  \ddots    & \ddots       & \\
          &          & d_{2,\varepsilon}          & 0       &            &            &      d_{1,\varepsilon}     &  1 \\
\end{array} \right]}.
\end{equation}
Then the problem is reduced to estimate the singular values of $\hat{\bb{L}}_{\varepsilon}$.

Step 2: As in the proof of Theorem \ref{thm:diffusive}, we first consider the simple case of $\varepsilon=0$. For simplicity, we omit the subscript 0 ($\varepsilon=0$) in the following. In this case,
\begin{align*}
& c_1 = \frac{1}{1+4\tilde{\beta} \sin^2 \frac{kh}{2} } \Big( -1-4\tilde{\beta} \sin^2 \frac{kh}{2} + ( 2\beta + \tilde{\beta} ) \sin^2 (kh) \Big),\\
& c_2 = d_1 = 0, \\
& d_2 = {\rm i} \frac{1}{1+4\tilde{\beta} \sin^2 \frac{kh}{2} } \Big( \beta  - 2\tilde{\beta}\sin^2 (kh) \Big)\sin (kh).
\end{align*}
When $\beta\le 1$, one has
\begin{align*}
|d_2|
& \le \frac{\beta + 2\tilde{\beta}\sin^2 (kh)}{1+4\tilde{\beta} \sin^2 \frac{kh}{2} }
 = \frac{\beta + 8\tilde{\beta}\sin^2 \frac{kh}{2}\cos^2 \frac{kh}{2}}{1+4\tilde{\beta} \sin^2 \frac{kh}{2} }
  \le \frac{1 + 8\tilde{\beta}\sin^2 \frac{kh}{2}}{1+4\tilde{\beta} \sin^2 \frac{kh}{2} } \le 2.
\end{align*}
In view of the inequality \eqref{sigmge},
\[\sigma_{\min}(\hat{\bb{L}}) \ge \sigma_{\min}(\hat{\bb{L}}_{11})/3.\]
Thus it suffices to bound $\sigma_{\min}(\hat{\bb{L}}_{11})$.

Step 3: By definition, $\sigma_{\min}(\hat{\bb{L}}_{11}) = 1/\sigma_{\max}(\hat{\bb{L}}_{11}^{-1})$, where
\[\hat{\bb{L}}_{11}^{-1} =
\begin{bmatrix}
1          &              &                 &          &            \\
-c_1        & 1       &                 &          &        \\
(-c_1)^2      &  \ddots      &  \ddots         &          &\\
   \vdots       &  \ddots      &  \ddots         &  \ddots  &   \\
(-c_1)^{N_t-1} &   \cdots     &   (-c_1)^2     &  -c_1 &  1      \\
\end{bmatrix}.\]
By Lemma \ref{lem:Qi},
\[\sigma_{\max}(\hat{\bb{L}}_{11}^{-1}) \le 1 + |c_1| + \cdots + |c_1|^{N_t-1} .\]
If $\beta \le 1/2$, then
\[
\frac{( 2\beta + \tilde{\beta} ) \sin^2 (kh)}{1+4\tilde{\beta} \sin^2 \frac{kh}{2} }
= \frac{2\beta \sin^2 (kh) + 4\tilde{\beta} \sin^2 \frac{kh}{2} \cos^2 \frac{kh}{2} }{1+4\tilde{\beta} \sin^2 \frac{kh}{2} }
\le 1,
\]
which yields
\[|c_1|
 = \left | 1 - \frac{( 2\beta + \tilde{\beta} ) \sin^2 (kh)}{1+4\tilde{\beta} \sin^2 \frac{kh}{2} } \right|
 \le 1\]
and hence
\[\sigma_{\max}(\hat{\bb{L}}_{11}^{-1}) \le N_t = 1/\tau .\]
Thus,
\[\sigma_{\min}(\hat{\bb{L}}) \ge \tau/3. \]

Step 4: Let $\bb{E} = \hat{\bb{L}}_{\varepsilon} - \hat{\bb{L}}$, where
\[ \bb{E} = {\scriptsize \left[\begin{array}{cccc:cccccccc}
0         &          &               &         &   0      &         &              &  \\
\tilde{c}_1      &      0   &               &         &     \tilde{c}_2       & 0     &              &   \\
          & \ddots   & \ddots        &         &            &    \ddots     &    \ddots    & \\
          &          & \tilde{c}_1           & 0       &            &         &      \tilde{c}_2        & 0\\
 \hdashline
0         &          &               &         & 0          &         &              & \\
\tilde{d}_2       & 0        &               &         &      \tilde{d}_1   & 0       &              & \\
          & \ddots   & \ddots        &         &            & \ddots  & \ddots       & \\
          &          & \tilde{d}_2           & 0       &            &          &      \tilde{d}_1    &  0 \\
\end{array} \right]},\]
with
\begin{align*}
& \tilde{c}_1 = c_{1,\varepsilon}-c_1 = \frac{\beta^2}{1+4\tilde{\beta} \sin^2 \frac{kh}{2} } \Big( \frac{\nu}{\gamma}-\frac{1}{\tau}  \Big) \sin^2 (kh),\\
& \tilde{c}_2 = c_{2,\varepsilon}-c_2 = {\rm i} \frac{\beta}{\gamma} \sin (kh), \\
& \tilde{d}_1 = d_{1,\varepsilon}-d_1 = \frac{2\beta}{\gamma} \sin^2 \frac{kh}{2}, \\
& \tilde{d}_2 = d_{2,\varepsilon}-c_2 = {\rm i} \frac{ 2\beta^2 }{1+4\tilde{\beta} \sin^2 \frac{kh}{2} } \Big(\frac{1}{\tau} - \frac{\nu}{\gamma} \Big) \sin^2 (kh)\sin (kh).
\end{align*}
Then,
\begin{align*}
\|\bb{E}\|
& \le |\tilde{c}_1| + |\tilde{c}_2| + |\tilde{d}_1| + |\tilde{d}_2|
\le 2 \Big(\frac{\varepsilon}{\tau+\varepsilon}\frac{1+\tau}{\tau} + \frac{\varepsilon}{\tau+\varepsilon} \Big)
= \frac{\varepsilon}{\tau+\varepsilon} \frac{2+4\tau}{\tau}.
\end{align*}

Step 5: According to the above calculations,
\begin{align*}
\sigma_{\max}(\hat{\bb{L}}_{\varepsilon})
& = \|\hat{\bb{L}}_{\varepsilon}\|
  \le 1 + |c_{1,\varepsilon}| + |c_{2,\varepsilon}| + |d_{1,\varepsilon}| + |d_{2,\varepsilon}| \\
& \le 1 + |\tilde{c}_1| + |\tilde{c}_2| + |\tilde{d}_1| + |\tilde{d}_2| + |c_1| + |c_2| + |d_1| + |d_2| \\
& \le \frac{\varepsilon}{\tau+\varepsilon} \frac{2+4\tau}{\tau} + 4.
\end{align*}
The desired estimates follow from the Weyl's inequality.

(2) Now we consider the time complexity. In what follows, we set $\tau = \mathcal{O}(h)$. The first step of the penalized diffusive relaxation scheme is to solve
\[
\begin{cases}
 \widetilde{\bb{B}}\bb{u}^* = \bb{u}^n + \tilde{\beta}\bb{b}^n, \\
 \gamma\bb{v}^* = \bb{v}^n - \nu \bb{A}\bb{u}^* - \nu \frac{\beta}{2} \widetilde{\bb{b}}^n.
\end{cases}
\]
The matrix $\widetilde{\bb{B}}$ is positive definite with order $n_u = N_x-1$. Let $\tilde{\beta} = \tau/h^2 \sim h^{-1}$. The condition number of $\widetilde{\bb{B}}$ is then given by
\[\kappa = \frac{1+\tilde{\beta}\pi^2}{1+8\tilde{\beta}h^2} = \mathcal{O}(h^{-1}) = \mathcal{O}(N_x). \]
The time complexity of obtaining $\bb{u}^*$ using the CG method is
\[\mathcal{O}(n_u \tilde{s} \sqrt{\kappa} \log (1/\delta)) = \mathcal{O}(N_x^{1.5} \log N_x),\]
where $\tilde{s}=3$ is the sparsity number of $\widetilde{\bb{B}}$. And the time complexity of obtaining $\bb{v}^*$ is $\mathcal{O}(N_x)$. Thus the time complexity of the first step is $\mathcal{O}(N_x^{1.5} \log N_x )$.

The second step is to solve
\[
\begin{cases}
\bb{u}^{n+1} = (\bb{B}+\widetilde{\bb{B}}-\bb{I}) \bb{u}^* - \bb{A}\bb{v}^* + \frac{\beta}{2}(\bb{b}^n - \widetilde{\bb{c}}^n)-\tilde{\beta}\bb{b}^n, \\
\bb{v}^{n+1} = \bb{B} \bb{v}^* - \bb{A}\bb{u}^* + \frac{\beta}{2}(\bb{c}^n - \widetilde{\bb{b}}^n).
\end{cases}
\]
It is easy to know that its time complexity is $\mathcal{O}(N_x)$. In summary, the time complexity of each iteration is $\mathcal{O}(N_x^{1.5}\log N_x)$,  thus the total running time is
\[C_{\text{penalized}} = \mathcal{O}(N_t N_x^{1.5}\log N_x ) = \mathcal{O}(N_x^{2.5}\log N_x),\]
where $\tau = \mathcal{O}(h)$ is used.

We now consider the quantum treatment. From the estimates of eigenvalues, the condition number $\kappa = \mathcal{O}(\tau^{-1}) = \mathcal{O}(N_t)$. Plugging them in \eqref{cpCAS}, one obtains
\begin{align*}
Q_{\text{penalized}}
 = \mathcal{O}(N_t \log (1/\delta) ) = \mathcal{O}(N_x \log N_x ) ,
\end{align*}
where $\tau = \mathcal{O}(h)$ is used in the last step.

The result for $\tau = \mathcal{O}(h^2)$ can be proved in similar way.
\end{proof}

As observed, all three quantum methods have an acceleration advantage when $\varepsilon \to 0$.

\subsection{Time complexity of the explicit scheme}

The explicit scheme is
\begin{equation}\label{explicitODE}
\begin{cases}
\frac{u_j^{n+1} - u_j^n}{\tau} + \frac{v_{j+1}^n-v_{j-1}^n}{2h} - \frac{h}{2\sqrt{\varepsilon}}\frac{u_{j-1}^n-2u_j^n+u_{j+1}^n}{h^2}= 0, \\
\frac{v_j^{n+1} - v_j^n}{\tau} + \frac{u_{j+1}^n-u_{j-1}^n}{2h \varepsilon} - \frac{h}{2\sqrt{\varepsilon}}\frac{v_{j-1}^n-2v_j^n+v_{j+1}^n}{h^2}
    = - \frac{1}{\varepsilon} v_j^n ,
\end{cases}
\end{equation}
which is written in vector as
\begin{equation*}
\begin{cases}
\frac{\bb{u}^{n+1} - \bb{u}^n}{\tau}  + \frac{1}{2h} \bb{M}_h \bb{v}^n - \frac{1}{2h \sqrt{\varepsilon}} \bb{L}_h \bb{u}^n - \frac{1}{2h} (\bb{b}^n /\sqrt{\varepsilon} - \widetilde{\bb{c}}^n) = \bb{0}, \\
\frac{\bb{v}^{n+1} - \bb{v}^n}{\tau}  + \frac{1}{2h\varepsilon} \bb{M}_h \bb{u}^n - \frac{1}{2h \sqrt{\varepsilon}} \bb{L}_h \bb{v}^n + \frac{1}{2h} (\widetilde{\bb{b}}^n/\varepsilon - \bb{c}^n/\sqrt{\varepsilon})
   =-\frac{1}{\varepsilon}\bb{v}^n .
\end{cases}
\end{equation*}
Let $\beta = \tau/h$. One has
\begin{equation} \label{explicithyperheat}
\begin{cases}
-\bb{B}\bb{u}^n + \bb{u}^{n+1} + \bb{A}\bb{v}^n = \bb{f}^{n+1}, \\
-(\bb{B}-\tau/\varepsilon\bb{I})\bb{v}^n +  \bb{v}^{n+1} + \varepsilon^{-1}\bb{A}\bb{u}^n = \bb{g}^{n+1},
\end{cases}
\end{equation}
where
\[\bb{A} = \frac{\beta}{2}\bb{M}_h, \qquad \bb{B} = \bb{I} + \frac{\beta}{2\sqrt{\varepsilon}} \bb{L_h},\]
\[\bb{f}^{n+1} = \frac{\beta}{2}(\bb{b}^n/\sqrt{\varepsilon}-\widetilde{\bb{c}}^n), \qquad \bb{g}^{n+1} = \frac{\beta}{2}(\bb{c}^n/\sqrt{\varepsilon} - \widetilde{\bb{b}}^n/\varepsilon). \]

Introduce the following notations
\begin{equation}\label{UVS}
\bb{U} = [\bb{u}^1; \cdots ;\bb{u}^{N_t}], \qquad \bb{V} = [\bb{v}^1; \cdots ;\bb{v}^{N_t}], \qquad \bb{S} = [\bb{U}; \bb{V}].
\end{equation}
Then one obtains the linear system
\begin{equation}\label{APexplicit}
\bb{L}_{\text{explicit}} \bb{S} = \bb{F}_{\text{explicit}},
\end{equation}
where $\bb{L}_{\text{explicit}} = (\bb{L}_{ij})_{2\times 2}$, $\bb{F}_{\text{explicit}} = [\bb{F}_1; \bb{F}_2]$, and
\[\bb{L}_{11} =
\begin{bmatrix}
\bb{I}  &            &           &            \\
-\bb{B} & \bb{I}     &           &            \\
        &\ddots      & \ddots    &    \\
        &            & -\bb{B}   & \bb{I}     \\
\end{bmatrix}
, \qquad
\bb{L}_{12} =
\begin{bmatrix}
\bb{O}  &            &           &            \\
\bb{A} & \bb{O}     &           &            \\
        &\ddots      & \ddots    &    \\
        &            & \bb{A}   & \bb{O}     \\
\end{bmatrix},\]
\[\bb{L}_{21} =
\begin{bmatrix}
\bb{O}  &            &           &            \\
\varepsilon^{-1}\bb{A} & \bb{O}     &           &            \\
        &\ddots      & \ddots    &    \\
        &            & \varepsilon^{-1}\bb{A}   & \bb{O}     \\
\end{bmatrix},\quad
\bb{L}_{22} =
\begin{bmatrix}
\bb{I}  &            &           &            \\
-(\bb{B}-\tau/\varepsilon\bb{I}) & \bb{I}     &           &            \\
        &\ddots      & \ddots    &    \\
        &            & -(\bb{B}-\tau/\varepsilon\bb{I})   & \bb{I}    \\
\end{bmatrix},\]
\[\bb{F}_1 =
\begin{bmatrix}
\bb{f}^1 + \bb{B}\bb{u}^0-\bb{A}\bb{v}^0 \\
\bb{f}^2  \\
\vdots\\
\bb{f}^{N_t} \\
\end{bmatrix}, \qquad
\bb{F}_2 =
\begin{bmatrix}
\bb{g}^1 - \varepsilon^{-1}\bb{A}\bb{u}^0 + (\bb{B}-\tau/\varepsilon\bb{I})\bb{v}^0\\
\bb{g}^2  \\
\vdots\\
\bb{g}^{N_t} \\
\end{bmatrix}.
\]

\begin{theorem}\label{thm:APexplicitcp}
Let $\tau = c_p\sqrt{\varepsilon} h$ and $h = \sqrt{\varepsilon}\delta$, where $1/\sqrt{111} \le c_p\le 2/(2+\delta)$ is a constant and $\delta \le 1$ is the error bound.
Denote $\bb{H}$ to be the Hermitian matrix corresponding to $\bb{L}_{\text{explicit}}$. Then the condition number and the sparsity of $\bb{H}$ satisfy
\[\kappa = \mathcal{O}( N_t ) = \mathcal{O}((\varepsilon\delta)^{-1})\qquad \mbox{and} \qquad s = \mathcal{O}(1).\]
The explicit scheme can be solved with time complexity
\[C_{\text{explicit}} = \mathcal{O}(\varepsilon^{-1.5}\delta^{-2}), \qquad   Q_{\text{explicit}} = \mathcal{O}( \varepsilon^{-1.5}\delta^{-1}\log (1/\delta) )\]
for the classical treatment and the quantum treatment, respectively.
\end{theorem}
\begin{proof}
Since the problem is linear, we can simply apply the discrete Fourier transform to characterize the singular values of the coefficient matrix $\bb{L}_{\text{explicit}}$ as done in the proof of Theorem \ref{thm:penalty}.

(1)  Considering the amplification factor $\varepsilon^{-1}$ in $\bb{L}_{21}$, we first reformulate the system \eqref{APexplicit} as
\[\begin{bmatrix} \bb{L}_{11}  & \sqrt{\varepsilon}^{-1} \bb{L}_{12} \\
\sqrt{\varepsilon} \bb{L}_{21}  &  \bb{L}_{22}
\end{bmatrix}\begin{bmatrix}
\sqrt{\varepsilon}^{-1}\bb{U}\\
\bb{V}
\end{bmatrix} = \begin{bmatrix}
\sqrt{\varepsilon}^{-1}\bb{F}_1\\
\bb{F}_2
\end{bmatrix},\]
with the coefficient matrix denoted by $\tilde{\bb{L}}_{\text{explicit}}$ in the following. This means we consider a linear system with new variables $\tilde{u} = u/\sqrt{\varepsilon}$ and $\tilde{v} = v$. Note that $\sqrt{\varepsilon}^{-1} \bb{L}_{12} = \sqrt{\varepsilon} \bb{L}_{21}$.

(2) Plugging the expressions of \eqref{fouriervariables} in \eqref{explicitODE}, one has
\[
\begin{cases}
\frac{\hat{u}^{n+1} - \hat{u}^n}{\tau} + \frac{({\rm e}^{{\rm i} kh}- {\rm e}^{-{\rm i}kh})\hat{v}^n}{2h} - \frac{h}{2\sqrt{\varepsilon}}\frac{({\rm e}^{-{\rm i} kh}-2 + {\rm e}^{{\rm i}kh})\hat{u}^n}{h^2} = 0, \\
\frac{\hat{v}^{n+1} - \hat{v}^n}{\tau} + \frac{({\rm e}^{{\rm i} kh}- {\rm e}^{-{\rm i}kh})\hat{u}^n}{2h\varepsilon} - \frac{h}{2\sqrt{\varepsilon}}\frac{({\rm e}^{-{\rm i} kh}-2 + {\rm e}^{{\rm i}kh})\hat{v}^n}{h^2} =
- \frac{1}{\varepsilon}\hat{v}^n ,
\end{cases}
\]
which can be written as
\[\begin{cases}
\hat{u}^{n+1} - c_1 \hat{u}^n - c_2 \hat{v}^n = 0\\
\hat{v}^{n+1} - d_1 \hat{v}^n - d_2 \hat{u}^n = 0
\end{cases}, \quad n = 0,1,\cdots, N_t-1\]
or
\begin{equation}\label{fourieroriginal}
\begin{bmatrix}\hat{u}^{n+1} \\ \hat{v}^{n+1} \end{bmatrix}
= \begin{bmatrix}c_1 & c_2 \\ d_2 & d_1 \end{bmatrix}
\begin{bmatrix}\hat{u}^n \\ \hat{v}^n \end{bmatrix} =: A \begin{bmatrix}\hat{u}^n \\ \hat{v}^n \end{bmatrix},
\end{equation}
where
\begin{align*}
 c_1 =  1 - \frac{2\beta}{\sqrt{\varepsilon}} \sin^2 \frac{kh}{2} , \quad
 c_2 = -{\rm i} \beta \sin (kh) , \quad
 d_1 =  1 - \frac{2\beta}{\sqrt{\varepsilon}} \sin^2 \frac{kh}{2} - \frac{\tau}{\varepsilon}, \quad
 d_2 = - {\rm i} \frac{\beta}{\varepsilon}\sin (kh).
\end{align*}
For the reformulated system, we should introduce new variables $\tilde{u} = \hat{u}/\sqrt{\varepsilon}$ and $\tilde{v} = \hat{v}$ and rewrite \eqref{fourieroriginal} as
\begin{equation}\label{uvfourierexplicit}
\begin{bmatrix}\tilde{u}^{n+1} \\ \tilde{v}^{n+1} \end{bmatrix}
= \begin{bmatrix}c_1 & c_2/\sqrt{\varepsilon} \\ \sqrt{\varepsilon} d_2 & d_1 \end{bmatrix}
\begin{bmatrix}\tilde{u}^n \\ \tilde{v}^n \end{bmatrix} = \begin{bmatrix}c_1 & \tilde{c}_2 \\ \tilde{d}_2 & d_1 \end{bmatrix}
\begin{bmatrix}\tilde{u}^n \\ \tilde{v}^n \end{bmatrix} =: \tilde{A} \begin{bmatrix}\hat{u}^n \\ \hat{v}^n \end{bmatrix},
\end{equation}
where $\tilde{c}_2 = \tilde{d}_2 = -\i \beta \sin(kh) /\sqrt{\varepsilon}$.
Let $a = |\sin(kh/2)| \in [0,1]$. One has
\[ \tilde{A} = \begin{bmatrix} 1-2\beta/\sqrt{\varepsilon} a^2 & \pm \i 2\beta/\sqrt{\varepsilon} a\sqrt{1-a^2} \\ \pm \i 2\beta/\sqrt{\varepsilon} a\sqrt{1-a^2} & 1-2\beta/\sqrt{\varepsilon} a^2 - \tau/\varepsilon\end{bmatrix}.\]

(3) The truncation error of the explicit scheme may be given by $\mathcal{O}(\tau + h/\sqrt{\varepsilon})$, so one can set $h = \sqrt{\varepsilon}\delta$, where $\delta$ is the error bound. Observing the term $\tau/\varepsilon$, we take $\tau = c_p \sqrt{\varepsilon} h$ or $\beta = \tau/h = c_p \sqrt{\varepsilon}$, where $c_p$ is a constant to be determined, hence
\[ \tilde{A} = \begin{bmatrix} 1-2 c_p a^2 & \pm \i 2 c_p a\sqrt{1-a^2} \\ \pm \i 2 c_p a\sqrt{1-a^2} & 1-2c_p a^2 - c_p\delta\end{bmatrix}.\]
One can verify that $\|\tilde{A}\| \le 1$ when $c_p$ and $\delta$ satisfy the given condition.

The argument is as follows. Let $\text{det}( \mu I - \tilde{A}\tilde{A}^\dag) = : \mu^2 - b \mu - c = 0 $ and denote $\mu_1$ and $\mu_2$ to be the two roots. A direct calculation gives
\[\tilde{A}\tilde{A}^\dag
= \begin{bmatrix} 4a^2 c_p^2 - 4a^2 c_p + 1 & -\i  \delta 2 a\sqrt{1-a^2} c_p^2 \\
-\i  \delta 2 a\sqrt{1-a^2} c_p^2 & (4a^2+4a^2\delta+\delta^2)c_p^2 - 2(2a^2+\delta)c_p + 1\end{bmatrix},\]
and
\[b = (8a^2 + 4a^2\delta + \delta^2) c_p^2  - (8a^2 + 2\delta) c_p + 2,\]
\begin{align*}
-c
& = 4(- a^4\delta^2 + 4a^4\delta + 4a^4 + 2a^2\delta^2) c_p^4 - 4(4a^4\delta+8a^4- a^2\delta^2 + 2a^2\delta) c_p^3\\
& \quad + (16a^4 + 12a^2\delta+ 8a^2 + \delta^2) c_p^2  - ( 8a^2  + 2\delta) c_p + 1.
\end{align*}
Recall the simple fact:
\[|\mu_i|\le 1~(i=1,2)  \qquad \Longleftrightarrow \qquad  |b|\le 1-c, \quad  |c|\le 1.\]
We claim that the right-hand side holds if $b$ and $-c$ are two monotonically decreasing functions on $[0,1]$ with respect to the variable $a$:
\begin{itemize}
\item If $-c$ is decreasing, then
\[-c \le -c(0) = (1-\delta c_p)^2 \le 1,\]
and
\begin{align*}
-c \ge -c(1)
& = (4c_p^4 + 4c_p^3 + c_p^2) \delta^2
 + (16 c_p^4-24c_p^3 + 12 c_p^2 - 2 c_p) \delta \\
& \quad  + 16c_p^4  - 32 c_p^3 + 24 c_p^2  -  8 c_p  + 1 := d_2 \delta^2 + d_1 \delta + d_0 =: f(\delta).
\end{align*}
  Noting that
  \[f_{\min}(\delta) = f\Big(-\frac{d_1}{2d_2}  \Big) = \frac{4d_0d_2 - d_1^2}{4d_2} = \frac{8c_p^3(2c_p - 1)^4}{d_2}  \ge 0,\]
  we obtain
 $-c  \ge -c(1) \ge 0$.
  That is, $|c| = -c \le 1$.
  \item If $b$ is decreasing, then
\[b \ge b(1) = c_p^2\delta^2 + (4-2c_p)\delta + 2(1-2c_p)^2 \ge 0,\]
\[b \le b(0) = 1 + (1-\delta c_p)^2.\]

  \item With the above estimates, one further gets
  \begin{align*}
  |b|+c
  & = b+c \le b(0) + c(1) \\
  & \le 1 + (1-\delta c_p)^2  - 16c_p^2(c_p^2-2c_p+8) - (4c_p-1)^2 \\
  & = 1 - (112-\delta^2) c_p^2 + 1 - 2\delta c_p -16 c_p^2(1-c_p)^2 -  (4c_p-1)^2.
  \end{align*}
  Obviously, $|b|+c\le 1$ if $- (112-\delta^2) c_p^2 + 1 \le 0$ or $c_p \ge 1/\sqrt{112-\delta^2} \ge 1/\sqrt{111}$.
\end{itemize}

It remains to verify that $b'(a)\le 0$ and $-c'(a)\le 0$. From $b'(a)\le 0$, one easily gets $c_p \le 2/(2+\delta)$.
A direct calculation gives
\begin{align*}
-c'(a)
& = 16(- a^3\delta^2+ 4a^3\delta + 4a^3 + a\delta^2) c_p^4
 + 8(- 8a^3\delta- 16a^3- a\delta^2- 2a\delta)c_p^3 \\
& \quad + 8(8a^3 + 3a\delta + 2a)c_p^2 - 16ac_p.
\end{align*}
From $-c'(a)\le 0$ one has
\begin{align}
2(-  c_p^2 \delta^2 - 4 c_p(1-c_p) \delta + 2c_p^2-8c_p+4)c_p  a^2 \nonumber\\
  \le 16  + (- 2 c_p^3 + c_p^2) \delta^2  + (2c_p^2 - 3c_p) \delta - 2c_p. \label{cpaa}
\end{align}
Let
\[f(\delta) = 2(-  c_p^2 \delta^2 - 4 c_p(1-c_p) \delta + 2c_p^2-8c_p+4)c_p  a^2,\]
\[g(\delta) = 16  + (- 2 c_p^3 + c_p^2) \delta^2  + (2c_p^2 - 3c_p) \delta - 2c_p.\]
One easily finds that $f'(\delta)\le 0$ and  $g'(\delta)\le 0$ when $\delta\le 1$ and $c_p\le 1$. Hence $f(\delta)\le f(0) \le 2c_p(c_p-2)\le 0$ and $g(\delta) \ge g(1) = 16 - 5c_p + 2c_p^3+c_p^2 \ge 0$. Therefore, \eqref{cpaa} is trivial when $\delta\le 1$ and $c_p\le 1$.

(4) We now bound the singular values of the reformulated matrix $\tilde{\bb{L}}_{\text{explicit}}$. The coefficient matrix $\tilde{\bb{L}}$ corresponding to \eqref{uvfourierexplicit} has the same form of \eqref{LFourier}. It is obvious that $\sigma_{\max}(\tilde{\bb{L}}) \lesssim 1$. For the smallest singular value, by definition, $\sigma_{\min}(\tilde{\bb{L}}) = 1/\|\tilde{\bb{L}}^{-1}\|$, so it suffices to give an upper bound of $\|\tilde{\bb{L}}^{-1}\|$. Noting that
\[\|\tilde{\bb{L}}^{-1}\| = \max_{\|\tilde{\bb{b}}\|\le 1} \|\tilde{\bb{L}}^{-1} \tilde{\bb{b}}\|, \qquad \tilde{\bb{b}} = [\tilde{\bb{f}}; \tilde{\bb{g}}],\]
where $\tilde{\bb{f}} = [\tilde{f}^1, \cdots, \tilde{f}^{N_t}]^T$ and $\tilde{\bb{g}} = [\tilde{g}^1, \cdots, \tilde{g}^{N_t}]^T$, one can determine the upper bound of $\|\tilde{\bb{L}}^{-1} \tilde{\bb{b}}\|$ by bounding the solution to the following linear system
\[\begin{bmatrix}\tilde{u}^{n+1} \\ \tilde{v}^{n+1} \end{bmatrix}
= \begin{bmatrix}c_1 & c_2/\sqrt{\varepsilon} \\ \sqrt{\varepsilon} d_2 & d_1 \end{bmatrix}
\begin{bmatrix}\tilde{u}^n \\ \tilde{v}^n \end{bmatrix} + \begin{bmatrix}\tilde{f}^n \\ \tilde{g}^n \end{bmatrix}
= \tilde{A}\begin{bmatrix}\tilde{u}^n \\ \tilde{v}^n \end{bmatrix} + \begin{bmatrix}\tilde{f}^n \\ \tilde{g}^n \end{bmatrix}.\]
Suppose the maximum is attained at $\tilde{\bb{b}}$. Let $\tilde{w}^n = [\tilde{u}^n, \tilde{v}^n]^T$ and $\tilde{b}^n = [\tilde{f}^n, \tilde{g}^n]^T$. One has $\tilde{w}^{n+1} = A \tilde{w}^n + \tilde{b}^n$, hence
\[\|\tilde{w}^n\| \le \|\tilde{A}\|^n \|\tilde{w}^0\| + \|\tilde{A}\|^{n-1} \|\tilde{b}^0\| +  \|\tilde{A}\|^{n-2} \|\tilde{b}^1\| + \cdots + \|\tilde{b}^{n-1}\|.\]
This gives $\|\tilde{w}^n\| \lesssim \|\tilde{w}^0\| + 1$ under the condition of the spatial and temporal steps since $\|\tilde{A}\| \le 1$ and $\|\tilde{\bb{b}}\|\le 1$. We then obtain
\[\|\tilde{\bb{L}}^{-1}\| = \|\tilde{\bb{L}}^{-1} \tilde{\bb{b}}\| = \|\tilde{\bb{w}}\|\lesssim N_t (\|\tilde{w}^0\| + 1),\]
where $\tilde{\bb{w}} = [\tilde{w}^1; \cdots; \tilde{w}^{N_t}]$, and thus $\sigma_{\min}(\tilde{\bb{L}}) \gtrsim 1/N_t$.

(5) From $\|\tilde{w}^n\| \lesssim \|\tilde{w}^0\| + 1$ we get
  \[ |\hat{u}^n| \le \sqrt{\varepsilon} \|\tilde{w}^n\| \lesssim \sqrt{\varepsilon} (\|\tilde{w}^0\| + 1) \le |\hat{u}^0| +  \sqrt{\varepsilon} (|\hat{v}^0| + 1),\]
which implies the stability for $u$. When performing the similar analysis below Eq.~\eqref{AP2Pre}, we may also deduce the stability for $v$, hence $\|A^n\| \le C$ for some constant $C$ independent of $\tau$ and $h$, where $A$ is given by \eqref{fourieroriginal}. Therefore, the analysis in the last step is also valid for the original linear system \eqref{APexplicit}, which implies $\sigma_{\min}(\bb{H}) \gtrsim 1/N_t$. It is obvious that   $\sigma_{\max}(\bb{H}) \lesssim 1/\sqrt{\varepsilon}$. We therefore obtain  $ \kappa(\bb{H}) \lesssim \mathcal{O}(N_t/\sqrt{\varepsilon}) = \mathcal{O}(\varepsilon^{-1.5}\delta^{-1})$.

(6) For the classical method, it is easy to get that the time complexity at each iteration is $\mathcal{O}(N_x)$, and hence the total run time is
\[C_{\text{explicit}} = \mathcal{O}(N_tN_x) = \mathcal{O}( (\varepsilon^{-1.5}\delta^{-2}).\]
For the quantum treatment, plugging the estimate of the condition number in \eqref{cpCAS}, one has
\begin{align*}
Q_{\text{explicit}}
 = \mathcal{O}(\kappa \log (1/\delta) ) = \mathcal{O}( \varepsilon^{-1.5}\delta^{-1} \log (1/\delta) ) .
\end{align*}
This completes the proof.
\end{proof}

As one can see, for both classical and quantum algorithms, the explicit methods have time complexity that depends on $\varepsilon$.

\section{Conclusions}

We studied the time complexities of finite difference methods for solving several time-dependent ODE and PDE problems, and in particular, a multiscale hyperbolic systems, in the setting of quantum computing,   The detailed results are summarized in Tab.~\ref{tab:summary} and Tab.~\ref{tab:summarymultiscale} for comparison.

There have been large strides in the development of QLSAs that show computational advantages over classical algorithms in certain regimes and under certain assumptions. To this end, it is natural to numerically solve the PDEs by combining the classical discretizations with QLSAs, especially for high-dimensional problems. This has been widely studied in recent years.

In this paper, we address issues in quantum algorithms for ODEs and high dimensional PDEs relevant for the concrete deployment of these algorithms. Firstly, we are interested in whether explicit--which is conditionally stable--and implcit--which are unconditionally stable--make any differences when one uses
quantum algorithms to solve these equations. For the heat equation we find that the different discretizations (forward Euler and Crank-Nicolson) don't make any difference when counting the cost of  the quantum algorithms.  We also present a unified analysis of the time complexity for the upwind discretization of the first order hyperbolic equation, which is a popular numerical scheme for such systems and the analysis may be generalized to other difference schemes, for example, the Lax-Friedrichs scheme and their high order approximations \cite{leveque2007finite}.

In addition to high-dimensional problems, the main focus of the paper is on quantum advantage of efficient multiscale methods for multiscale physical problems, which is important for applications  not yet studied in the quantum computing community.
We discuss in detail of the quantum difference methods for solving a prototype  multiscale problem and our results show that the Asymptotic-Preserving schemes, a popular multiscale framework for multiscale problems \cite{JinActa} -- are equally important in quantum computing since they allow the computational costs for quantum algorithms to be {\it independent} of the small physical scaling parameters. This also suggests that one should take full advantage of state-of-the-art multiscale classical algorithms when designing quantum algorithms for multiscale PDEs.

The approach and analysis in this paper can be extended to other more complex problems, such as the multiscale time-dependent transport equations in \cite{Jin-Pareschi-Toscani-2000}.

\begin{table}
  \begin{threeparttable}
  \centering
  \caption{Time complexities of classical and quantum difference methods for the simple ODE and the high-dimensional PDEs} \label{tab:summary}
  \setlength{\tabcolsep}{1.0mm}{
\begin{tabular}[\textwidth]{|c|cc|ccc|}
\hline
\textbf{Equation}   & \multicolumn{2}{c|}{\textbf{Classical difference methods}}  & \multicolumn{3}{c|}{\textbf{Quantum difference methods }} \\
\hline
  \multirow{2}{*}{ODE}
   & Forward Euler & C-N
   &   \multicolumn{2}{c}{Forward Euler}   & C-N            \\
   & $\mathcal{O}(N_t^2)$   & $\mathcal{O}(N_t)$
   &  \multicolumn{2}{c}{$ \mathcal{O}(N_t^2\log N_t)$}   & $\mathcal{O}(N_t\log N_t)$     \\
   \hline
  \multirow{2}{*}{Heat equation}
   & Explicit       & C-N
   &  \multicolumn{2}{c}{Explicit}        & C-N   \\
   & $\mathcal{O}(d^2 N_x^{d+2})$   & $\mathcal{O}(d^2N_x^{d+0.5})$
   &  \multicolumn{2}{c}{$ \mathcal{O} ( d N_x^2 \log(N_x^2/d)  )$}   & $\mathcal{O} ( d N_x^2 \log(N_x^2/d)  )$     \\
   \hline
    Hyperbolic equation
   & \multicolumn{2}{c|}{$\mathcal{O}(d^2N_x^{d+1})$}
   & \multicolumn{3}{c|}{$\mathcal{O}( d^2 N_x \log( N_x/d) ) $}   \\
 \hline
\end{tabular}}
  \end{threeparttable}
\end{table}

\begin{table}
  \begin{threeparttable}
  \centering
   \caption{Time complexities of classical and quantum difference methods for the multiscale problem} \label{tab:summarymultiscale}
  \setlength{\tabcolsep}{1.0mm}{
\begin{tabular}[\textwidth]{|c|cc|ccc|}
\hline
\textbf{Equation}   & \multicolumn{2}{c|}{\textbf{Classical difference methods}}  & \multicolumn{3}{c|}{\textbf{Quantum difference methods }} \\
\hline
  \multirow{5}{*}{\begin{tabular}[c]{@{}c@{}}Multiscale telegraph\\ equation\end{tabular}}
   & IMEX  $(\tau \sim h^2 )$      & $\mathcal{O}(\delta^{-3}) $
   & \multicolumn{3}{c|}{$\mathcal{O}(\delta^{-2}\log \delta^{-1})$}    \\
  \cline{2-6}
   & Relaxation $(\tau \sim h^2)$            & $\mathcal{O}(\delta^{-3}) $
   & \multicolumn{3}{c|}{$\mathcal{O}(\delta^{-2}\log \delta^{-1} )$}    \\
  \cline{2-6}
   & Penalized $(\tau \sim h^2)$         & $\mathcal{O}(\delta^{-3}\log \delta^{-1}) $
   & \multicolumn{3}{c|}{$\mathcal{O}(\delta^{-2}\log \delta^{-1} )$}    \\
  % \cline{2-6}
   & Penalized $(\tau \sim h)$        & $\mathcal{O}(\delta^{-2.5}\log \delta^{-1}) $
   & \multicolumn{3}{c|}{$\mathcal{O}(\delta^{-1} \log \delta^{-1})$}     \\
  \cline{2-6}
   & Explicit   $(\tau \sim \sqrt{\varepsilon}h)$      & $\mathcal{O}(\varepsilon^{-1.5}\delta^{-2}) $
   & \multicolumn{3}{c|}{$\mathcal{O}( \varepsilon^{-1.5}\delta^{-1} \log \delta^{-1} )$}\\
 \hline
\end{tabular}}
  \begin{tablenotes}
    \item[] $\varepsilon$ is the relaxation time or the scaling parameter, and $\delta$ is the error bound.
  \end{tablenotes}
  \end{threeparttable}
\end{table}

\subsection*{Declaration of competing interest}

The authors declare that they have no known competing financial interests or personal relationships that could have
appeared to influence the work reported in this paper.

\section*{Acknowledgement}
The authors wish to express their sincere appreciation to the referees for valuable comments and suggestions which
greatly improved an earlier version of the paper.
SJ was partially supported by the NSFC grant No.~12031013 and the Shanghai Municipal Science and Technology Major Project (2021SHZDZX0102).  NL acknowledges funding from the NSFC International Young Scientists Project (no.~12050410230), the Science and Technology Program of Shanghai (no.~21JC1402900) and the Natural Science Foundation of Shanghai grant 21ZR1431000. YY was partially supported by China Postdoctoral Science Foundation (no. 2022M712080).

%\end{CJK*}

\end{document}